\documentclass[regno,11pt]{amsart}
\usepackage{mymacros}

\usepackage{amsmath}
\usepackage{amssymb}

\usepackage{color}
\usepackage[margin=3.5cm]{geometry}
\usepackage{tikz}
\usepackage[hang,small,bf]{caption}
\usepackage[subrefformat=parens]{subcaption}

\title{An introduction to sheaf quantization}
\author{Tatsuki Kuwagaki}
\dedicatory{Dedicated to Professor Yoshitsugu Takei's Kanreki}

\begin{document}
%
% The text goes here.  
% Be sure to use the appropriate "theorem-like" environment as 
% is the following examples.  Never use plain TeX commands for these, as
% they will cause interference with the styles of other papers. 

\maketitle

%\tableofcontents      %optional
\begin{abstract}      %optional
The notion of sheaf quantization has many faces: an enhancement of the notion of constructible sheaves, the Betti counterpart of Fukaya--Floer theory, a topological realization of WKB-states in geometric quantization. 
The purpose of this note is to give an introduction to the subject. 
\end{abstract}

\section{Introduction}
Mikio Sato, one of the founders of microlocal analysis, emphasized the importance of considering the ``microlocal" (i.e., cotangent) direction, and the subject~\cite{SKK} was formulated in the language of symplectic geometry. The Riemann-Hilbert correspondence ~\cite{KashiwaraRH}~\footnote{A related foundation is given in \cite{KashiwaraKawai}, another proof is given by Mebkhout~\cite{Mebkhout}.} gave a dictionary linking analysis and the topological object ``sheaves". Kashiwara-Schapira~\cite{Microlocalstudy, KashiwaraSchapira} developed microlocal sheaf theory, which looks at sheaves in the cotangent directions and can be considered as a sheaf-theoretic counterpart of microlocal analysis under the RH-correspondence. Although the RH correspondence is valid only in the analytic framework, microlocal sheaf theory also works in a broader framework, in particular for real manifolds.

Kashiwara-Schapira's paper~\cite{Microlocalstudy} was published in 1985. Around the same time, symplectic topology gradually began to gain a foothold in mathematics, including Gromov's pseudo-holomorphic curves~\cite{Gromov}, Floer's approach to the Arnold conjecture~\cite{Floer}, and Fukaya category~\cite{Fukayacat}. A priori, their study of symplectic topology seems far removed from microlocal analysis because of the very different methods of investigation. I imagine that few, if any, foresaw that microlocal sheaf theory had to have anything to do with symplectic topology.\footnote{Some key pieces to envision the connection slowly came together. For example, the generalization of the Dubson-Kashiwara index theorem to Ext of constructible sheaves is shown to be related with Lagrangian intersections~\cite{GrinbergMacpherson}. Lagrangian intersection theory in symplectic topology is related to Fukaya category. There were several people who imagined a connection between Fukaya category and deformation quantization (Feigin, Bresler-Soibelman~\cite{BresslerSoibelman}).}

In such a situation, the works of Tamarkin~\cite{Tam} (apperaed on arXiv in 2008) and Nadler--Zaslow~\cite{NZ} (apperaed on arXiv in 2006) were a major breakthrough. Both papers suggest that microlocal sheaf theory is useful for symplectic topology, or, more naively, that the relationship between microlocal sheaf theory and symplectic geometry is very deep.

After their works, the concept of sheaf quantization was born. At first, it appeared as a tool to study non-conic Lagrangian submanifolds sheaf-theoretically. And then, it gradually becomes clear that sheaf quantization is a very rich concept and has many facets: (1) An enhancement of the notion of constructible (or monodromy) sheaves~\cite[\dots]{DK, kuwagaki2020sheaf, kuwagakihRH}, (2) the Betti counterpart of Fukaya--Floer theory~\cite[\dots]{Tam, GKS, Guillermou, JT, AsanoIkepersitence} etc.\ , (3) topological realization of WKB-states~\cite[\dots]{Tam, kuwagaki2020sheaf}.

In this note, I would like to introduce the concept of sheaf quantizations. In \S 2, I will explain the concept of sheaf quantization in a very vague way. Since the concept of sheaf quantization is still under development, it would be helpful to know the general goals. In \S 3, I will introduce a primitive version of the idea ``sheaf quantization version 0". This is just another name for a very well-known concept, (weakly) constructible sheaves. Following Nadler--Zaslow's work, we can view it in the symplectico-topological context. The goal of this section is to give an idea connecting sheaves to symplectic topology. In \S 4, the notion of sheaf quantization is introduced and I will explain several ideas behind its definition and applications. In this paper, I call this common sheaf quantization ``sheaf quantization version 1", to distinguish it from the version introduced in the following section. In \S 5, I will introduce sheaf quantization version 2, which is necessary to go beyond exact symplectic geometry. In \S 6 (Appendix), I also explain some key concepts around this study. If the reader does not know anything about microlocal sheaf theory, symplectic geometry, and wrapped Fukaya category, it is recommended to start with this section.

The followings are disclaimers: First, due to the structure of the paper, the topics are not presented in chronological order. Second, proofs are largely omitted, and some definitions and statements may not be stated accurately. If you want to know more, you can consult lecture notes/surveys that explain the relationship between sheaves and symplectic geometry include Viterbo~\cite{viterbo}, Guillermou--Schapira~\cite{GuiS}, and original literatures. Relevant papers will be introduced in the text as appropriate. Third, although I will mention several applications of the ideas, there are many other interesting applications of the sheaf-theoretic interpretation of Fukaya category, which are not introduced in this paper: Legendrian invariants~\cite{STZ}, Cluster varieties~\cite{STWZ, CasalsZaslow}, Lagrangian fillings~\cite{STW, GSW, CasalsGao}, Langlands program~\cite{BenZviNadler}, and so on. Of course, this list of literatures is far from complete.

\subsection*{Acknowledgment}
My knowledge of microlocal sheaf theory and symplectic geometry are largely due to many years of discussions with Yuichi Ike, Takahiro Saito, Fumihiko Sanda, and Vivek Shende. I would like to thank them here. Most of the materials presented here also consist of other things I have learned from various people. I am sorry that I cannot thank them all here.

Some parts of this note were prepared for the talks at (virtual) Yale university and (virtual) Academia Sinica. I'd like to thank the organizers of the seminars for their invitations. I thank Yuichi Ike, Pierre Schapira, Bingyu Zhang, and the anonymous referee for their comments on the early version.

This paper will be submitted to the proceedings of Takei 60, and I would also like to thank Takei-sensei for always kindly listening to me, an outsider in WKB analysis. I would also like to thank the organizers of the conference for inviting me.

\section{Sheaf quantization vaguely}
\subsection{Fukaya category}
One quick way to understand what sheaf quantization is is through the path of first learning about Fukaya category\footnote{As we will see later, this is only one of the motivations.}. It is not possible to go into the details of the definition of Fukaya category, but here is a rough outline. For relevant details, see, e.g., \cite{Auroux, FOOO, Seidel, Sectorialdescent}.

Let $(X, \omega)$ be a symplectic manifold. The Fukaya category of $(X, \omega)$ is a category (or, more precisely, a category with an appropriate enhancement, such as an $A_\infty$-category) described as follows:.
\begin{itemize}
    \item Object: Each object is a Lagrangian brane of $X$. It is determined by a Lagrangian submanifold of $X$ with brane data on it. If the reader does not know what brane data is, please ignore it this time. Anyway, it is an additional structure on a Lagrangian.
    \item Morphism: With a choice of the base ring $\bK$, the Fukaya category is defined as a $\bK$-linear category. The space of morphisms of Lagrangian branes from $L_1$ to $L_2$ is a free module generated by the intersections of their underlying Lagrangian submanifolds.
    \item Composition of morphisms: Suppose there are Lagrangian branes named $L_1, L_2, L_3$. Consider $p_1\in L_1\cap L_2$ and $p_2\in L_2\cap L_3$ as $p_1\in \Hom(L_1, L_2)$ and $p_2\in \Hom(L_2, L_3)$. Then the composition is defined as
    \begin{equation}\label{eqn:2.1}
    p_2\circ p_1=\sum_{p_3\in L_1\cap L_3}\#\cM(L_1, L_2, L_3, p_1, p_2,p_3)\cdot p_3.
    \end{equation}
    Let us now roughly describe what $\cM(L_1, L_2, L_3, p_1, p_2,p_3)$ is. First, we take one almost complex structure which makes $\omega$ into an almost K\"ahler form. Take three points $a_1, a_2, a_3$ on the boundary of the unit closed disc $D$ in $\bC$. Then $\cM(L_1, L_2, L_3, p_1, p_2,p_3)$ is defined as the moduli space of pseudo-holomorphic maps from $D$ to $M$ (i.e. consistent with the respective almost complex structures) that sends $a_i$ to $p_i$ and whose line segment in $\partial D$ connecting $a_i$ to $a_{i+1}$ lies $L_i$. More precisely, there is a restriction to the Maslov degree, which we omit here. Although the definition depends on the chosen almost complex structure, it is known that the quasi-equivalence class of the resulting $A_\infty$-category does not depend on it.
    \item 0-th and higher-order compositions of morphisms: Defining compositions of morphisms as above does not satisfy the associative laws. Therefore, for two, four or more objects and intersections, definitions similar to the above can be used to define a differential (0-th composition) and higher-order compositions.
    Then, a structure arises to compensate for the breaking of the associative law, which is called the $A_\infty$ structure.
\end{itemize}
In practice, there are many variants of the definition of Fukaya category. In particular, when the symplectic manifold is non-compact and exact, the objects are often restricted to exact Lagrangian branes (see \S \ref{wrapped} below). Also, we have to control perturbation of Lagrangians in the non-compact direction (cf. symplectic stops~\cite{Sylvan, Liouvillesector}). So, when we speak of a Fukaya category without specifying anything, we should generically imagine one that roughly satisfies the above situation.

\subsection{Sheaf quantization vaguely}
Sheaf quantization also has several variants as we will see below, but in general, the category formed by sheaf quantizations is expected to be equivalent to some kind of Fukaya category. It is worth pointing out that such equivalences between sheaf quantizations and Fukaya category are not dualistic equivalences, as opposed to homological mirror symmetry~\cite{hms}. Rather, just different representations of the same thing (i.e., A-model). More specifically, it is a generalization (or categorification) of the isomorphism between Morse homology and singular homology. In order to clarify the situation a little more in this perspective, the following terms are introduced.
\begin{dfn}
Let $(X, \omega)$ be a symplectic manifold. Let $\frakL(X, \omega)$ be the set of Lagrangian submanifolds of $X$ (or the set of Lagrangian branes in the situation where a Maslov structure is specified for $X$). Let $\cC$ be a category (or some enhanced category).
Then $\cC$ is a Lagrangian category of $(X, \omega)$ if a map $\frakL$ from the set of objects $\mathrm{Ob}(\cC)$ to $\frakL(X, \omega)$ is specified.
\end{dfn}
A Fukaya category is a Lagrangian category. 
Let us define the equivalence between Lagrangian categories  as follows.
\begin{dfn}
Let $\cC_1$ and $\cC_2$ be Lagrangian categories of $(X, \omega)$. Let $\frakL_i\colon \mathrm{Ob(\cC_i)}\rightarrow \frakL(X, \omega)$ be the structural morphisms. A category equivalence $F\colon \cC_1\rightarrow \cC_2$ is a Lagrangian category equivalence if there exists a set of generators $\lc L_a\rc_{a\in A}$ of $\cC_1$ such that $\frakL_2(F(L_a))=\frakL_1(L_a)$ for any $a\in A$.
\end{dfn}
This is also an ad hoc definition, and is used below in a suitably relaxed sense.

Of course, if we take the definitions literally, it is not helpful to distinguish non-dualistic equivalences from dualistic equivalences, since a category equivalent to a Lagrangian category can be a Lagrangian category. However, in the situation we are interested in, $\frakL$ is a priori defined without referring to any category equivalence to a Lagrangian category. Restricted to such categories, our terminology is a little helpful.

Now we define a slogan of sheaf quantization as a ``pseudo-definition".
We will define several categories of sheaf quantization starting in the next chapter, all of which satisfy the following properties, at least conjecturally. 
\begin{pdfn}[Sheaf quantization, vaguely]
Let $(X, \omega)$ be a symplectic manifold. Let $\cC$ be a Lagrangian category of $X$. We say $\cC$ is a category of sheaf quantizations if 
\begin{enumerate}
\item $\cC$ is a category of sheaves on the space associated to $X$,
    \item $\frakL$ is cooked out of the microsupport (Definition~\ref{microsupport}),
    \item $\cC$ is equivalent to a Fukaya category as a Lagrangian category of $(X, \omega)$.
\end{enumerate}
\end{pdfn}

\begin{rem}
A Fukaya category is a category defined by construction. It is experts' dream to characterize Fukaya category with something like the Eilenberg--Maclane axiom (there are some attempts, e.g., by Ganatra--Pardon--Shende). Such a characterization may give a way to rigidify the above pseudo-definition.
\end{rem}
\begin{rem}
As we will see later, $\frakL(X, \omega)$ is not enough. We eventually need to include singular Lagrangian and coisotropic submanifolds/branes.
\end{rem}

\section{Sheaf quantization version 0}
\subsection{Conic Lagrangians}
This version is usually not called sheaf quantization, but it is reasonable to think it as version 0 of sheaf quantization.

Let $M$ be a real manifold.\footnote{Sometimes, we need the real analyticity. See a relevant remark in Appendix.} As a symplectic manifold $(X, \omega)$, we consider the cotangent bundle $T^*M$ of $M$ with the standard symplectic structure. Let $\frakL(X, \omega)$ be the set of Lagrangians. We also allow singular Lagrangians here. Let us fix a field $\bK$ and our sheaves in this paper are sheaves valued in $\bK$-vector spaces.
\begin{dfn}
The category of sheaf quantization version 0 (or SQv0 in short) of $(T^*M,\omega)$ is the derived category $\Sh^{wc}(M)$ of weakly constructible sheaves (\S 6.1) on $M$.\footnote{This is only a tentative definition. Other options include: (1) the category of all the sheaves $\Sh(M)$, (2) the category of constructible sheaves $\Sh^c(M)$ (i.e., finite rank objects).}
\end{dfn}
For a weakly constructible sheaf $\cE$, its microsupport $\SS(\cE)$ defines a Lagrangian subset. That is, $\SS$ gives to $\Sh^{wc}(M)$ the structure of a Lagrangian category.
\begin{dfn}
We say an object $\cE$ of $\Sh^{wc}(M)$ is an SQv0 of $\SS(\cE)$ if $\SS(\cE)$ is a $C^0$-submanifold\footnote{It seems that if $\SS(\cE)$ is highly singular, the corresponding Lagrangian brane should not be considered as a brane on the singular Lagrangian, rather, a brane on the completion of the smooth Lagrangians. See Example~\ref{4.12} for related remarks.} of $T^*M$.
\end{dfn}

There is an important class of subcategories of $\Sh^{wc}(M)$. A subset of $T^*M$ is said to be conic if it is stable under the scaling of fibers of $T^*M$. Microsupport is always conic by its definition.
\begin{dfn}
Fix one conic subset $\bL$ in $T^*M$. Write $\Sh^{wc}_{\bL}(M)$ for the full subcategory spanned by the object of $\Sh^{wc}(M)$ whose microsupport is in $\bL$. Furthermore, we write $\Sh^w_{\bL}(M)$ for the full subcategory consisting of the compact~\footnote{The associated Yoneda module is cocontinuous.} objects of this category. An object of $\Sh^w_\bL(M)$ is called a wrapped constructible sheaf ~\cite{Nadlerwrapped}.
\end{dfn}

\begin{rem}
It is not so straight-forward to understand the meaning of taking compact objects. Other approaches are by localization ~\cite{Nadlerwrapped, IK, Sectorialdescent} or by wrapping ~\cite{Kuo}.
\end{rem}

\begin{exa}[Some examples of SQv0]
\begin{enumerate}
    \item For a point $x\in M$, let $T^*_xM$ be the cotangent fiber over $x$. An SQv0 of $T^*_xM$ is the skyscraper sheaf $\bK_x$.
    \item Let $T^*_MM$ be the zero-section of $M$. As an SQv0 of $T^*_MM$, we can take a local system on $M$.
    \item Let $i\colon N\hookrightarrow M$ be a submanifold. Then the conormal bundle $T^*_NM$ defines a Lagrangian submanifold of $T^*M$. The extrusion $i_*\cE$ of the local system $\cE$ on $N$ defines an SQv0 of $T^*_NM$.
\end{enumerate}
\end{exa}

As a very explicit relation between sheaf quantization and symplectic topology, the following theorem is of a fundamental importance:

\begin{thm}
[Nadler--Zaslow~\cite{NZ}, Ganatra--Pardon--Shende~\cite{MicrolocalMorse}]\label{theorem:NZGPS}
Let $S^*M$ be the contact boundary of the standard Liouville structure of $T^*M$. Determine $\partial \bL:=\bL\cap S^*M$ and regard it as a symplectic stop. In this case, there is a category equivalence
\begin{equation}
    \Sh^w_\bL(M)^{op}\cong \cW(T^*M, \partial\bL)
\end{equation}
where $op$ is the opposite category and $\cW(T^*M, \partial\bL)$ is the wrapped Fukaya category ~\cite{Sylvan, Sectorialdescent} of $(T^*M, \partial\bL)$, a kind of Fukaya category (see \S \ref{wrapped}).
\end{thm}
The theorem in this form is proved by \cite{MicrolocalMorse}. There is a closely related earlier version by Nadler--Zaslow~\cite{NZ}, which has been very influential for decades. More on precedents can be found in the reference of loc.\ cit. This category equivalence is a Lagrangian category equivalence in a certain sense. Let us see it in some examples.

\begin{exa}\label{exampleofsheafLagcorres1}
Let us consider the case when $M=S^1=\bR/\bZ$. The open intervals $(0,1/2)$ and $(1/2, 1)$ in $\bR$ are homeomorphically projected down to $S^1$. We denote the zero-extension of the constant sheaves over the intervals to $S^1$ by $\bK_{(0,1/2)}$ and $\bK_{(1/2, 1)}$ respectively.

We set $\bL:=T^*_0S^1\cup T^*_{1/2}S^1\cup T^*_{S^1}S^1$. Then the category $\Sh^w_{\bL}(S^1)$ is generated by skyscraper sheaves $\bK_0, \bK_{1/2}$ and the aforementioned sheaves $\bK_{(0,1/2)}, \bK_{(1/2, 1)}$. Under the above equivalence, the corresponding Lagrangians are $T^*_0 S^1, T^*_{1/2}S^1, T^{*,-}_{0}S^1\cup (0,1/2)\cup T^{*,+}_{1/2}S^1, T^{*,-}_{1/2}S^1\cup (1/2,1)\cup T^{*,+}_{1}S^1$, respectively, where $\pm$ indicates the positive/negative half of cotangent fibers. These are precisely microsupport of the corresponding sheaves.\footnote{Although the last two Lagrangians are singular, the singularities are very mild: One can easily set-up Floer theory for these Lagrangians. Or, one can choose obvious smoothings.}

\begin{figure}[htbp]
  \begin{minipage}[b]{0.45\linewidth}
    \centering
\begin{tikzpicture}
\begin{scope}[gray]
\begin{scope}[dashed]
\draw(2.4,-1.5) arc (0:180:1.2cm and 0.5cm);
\end{scope}
\draw(1.2,1.5) circle (1.2cm and 0.5cm);
\draw(0,-1.5) arc (180:360:1.2cm and 0.5cm);
\draw[-] (0,-1.5)-- (0,1.5);
\draw[-] (2.4,-1.5)-- (2.4,1.5);
\end{scope}
\draw[-,thick] (0,-1.5)-- (0,1.5);
\draw[-,thick] (2.4,-1.5)-- (2.4,1.5);
\draw[thick](1.2,0) circle (1.2cm and 0.5cm);
\draw (0,0.5) node[left] {$T^*_0S^1$};
\draw (2.4,0.5) node[right] {$T^*_{1/2}S^1$};
\draw (1.2,-0.5) node[below] {$T^*_{S^1}S^1$};
\end{tikzpicture}
\subcaption{$\bL$}
  \end{minipage}
  \begin{minipage}[b]{0.45\linewidth}
      \centering
\begin{tikzpicture}
\begin{scope}[gray]
\begin{scope}[dashed]
\draw(2.4,-1.5) arc (0:180:1.2cm and 0.5cm);
\end{scope}
\draw(1.2,1.5) circle (1.2cm and 0.5cm);
\draw(0,-1.5) arc (180:360:1.2cm and 0.5cm);
\draw[-] (0,-1.5)-- (0,1.5);
\draw[-] (2.4,-1.5)-- (2.4,1.5);
\end{scope}
\draw[-,thick] (0,-1.5)-- (0,0);
\draw[-,thick] (2.4,0)-- (2.4,1.5);
\draw[thick](0,0) arc (180:360:1.2cm and 0.5cm);
\draw (0,-0.5) node[left] {$T^{*,-}_0S^1$};
\draw (2.4,0.5) node[right] {$T^{*,+}_{1/2}S^1$};
\draw (1.2,-0.5) node[below] {$(0,1/2)$};
\end{tikzpicture}
\subcaption{$T^{*,-}_{0}S^1\cup (0,1/2)\cup T^{*,+}_{1/2}S^1$}
  \end{minipage}
  \caption{Example~\ref{exampleofsheafLagcorres1}}
\end{figure}
\end{exa}

\begin{exa}\label{exampleofsheafLagcorres2}
Let us consider a similar situation, but now we set $\bL:=T^*_0S^1\cup T^*_{S^1}S^1$. Then the category $\Sh^w_{\bL}(S^1)$ is generated by $\bK_0, \bK_{(0,1)}$. This time, the microsupport of $\bK_{(0,1)}$ is singular, which has a 4-valent point. In fact, $\bK_{(0,1)}$ has its microsupport as in Figure~\ref{figure:3.2}.

To have a Lagrangian category equivalence, we have to include the singular Lagrangian Figure~\ref{figure:3.2} in the Fukaya category. To do so, one can try to consider smoothing. But, reasonable smoothings of a 4-valent point are not unique. See Figure~\ref{fig:La} and Figure~\ref{fig:La'}. 

By studying the Nadler--Zaslow equivalence, one can conclude that $\bK_{(0,1)}$ corresponds to one of them. The sheaf corresponding to the other is $\bK_{[0,1]}$, which is push-forward of the constant sheaf $\bK_{[0,1]}$ to $\bR$ through the quotient morphism $[0,1]\rightarrow \bR\rightarrow \bR/\bZ$. The microsupport of $\bK_{[0,1]}$ is again Figure~\ref{figure:3.2}. 

From this observation, what we can learn is that we should enhance $\cL(X, \omega)$ so that each of its elements encode its smoothing. A clearer explanation can be found in Example~\ref{4.12}. Here we take another approach.

Instead, we consider the category $\Sh^w_{\bL}(S^1)$ as a quotient of the category of Example~\ref{exampleofsheafLagcorres1} by the subcategory generated by $\bK_{(0,1/2]}$ and $\bK_{[1/2,1)}$.\footnote{These objects are called microlocal skyscraper sheaves~\cite{Nadlerwrapped}. The corresponding objects on the Fukaya side are called linking disks The quotient description is known as ``stop removal".} For an object $\cE\in \Sh^w_{\bL}(S^1)$, instead of taking microsupport naively, let us first choose a representative and then take microsupport. Specifically, take $\bK_{(0, 1/2)}$ as a representative of $\bK_{(0,1)}$. This gives a structure of a Lagrangian category and the desired Lagrangian category equivalence.

\begin{figure}
    \centering
\begin{tikzpicture}
\begin{scope}[gray]
\begin{scope}[dashed]
\draw(2.4,-1.5) arc (0:180:1.2cm and 0.5cm);
\end{scope}
\draw(1.2,1.5) circle (1.2cm and 0.5cm);
\draw(0,-1.5) arc (180:360:1.2cm and 0.5cm);
\draw[-] (0,-1.5)-- (0,1.5);
\draw[-] (2.4,-1.5)-- (2.4,1.5);
\end{scope}
\draw[-,thick] (0,-1.5)-- (0,1.5);
\draw[thick](1.2,0) circle (1.2cm and 0.5cm);
\draw (1.2,-0.5) node[below] {$T^*_{S^1}S^1$};
\draw (0,0.5) node[left] {$T^*_0S^1$};
\end{tikzpicture}
\caption{Example~\ref{exampleofsheafLagcorres2}: $\bL=\SS(\bK_{(0,1)})=\SS(\bK_{[0,1]}$)\label{figure:3.2}}
\end{figure}
\end{exa}

\subsection{Why is it called ``quantization"?}
The reason why sheaf quantization version 0 is called ``quantization" is that it is a topological version (Betti version) of quantization.

The sheaf of rings of differential operators $\cD$ on a manifold $M$ is a specialization of the deformation quantization (see \S \ref{section:quantization}) of its cotangent bundle to $\hbar=1$. Also, a $\cD$-module can be considered as a deformation quantization of its characteristic variety in this sense. In particular, a holonomic $\cD$-module is a deformation quantization of a Lagrangian.

A holonomic $\cD$-module is said to be regular when its singularity is mild~(see, for example, \cite{HTT} for the definition). For this notion, the following classical RH correspondence (\cite{KashiwaraRH, KashiwaraKawai, Mebkhout}) is known.
\begin{thm}[Regular Riemann--Hilbert correspondence]
There is a category equivalence
\begin{equation}
    D^b_{rh}(\cD_M)\cong \Sh^{\bC c}(M),
\end{equation}
where the left-hand side is the derived category of regular holonomic $\cD$-modules and the right hand side is the derived category of $\bC$-constructible sheaves.
\end{thm}
Recently, it has been extended to the case of holonomic modules. We will later give a comment on this generalization.

Under the regular RH correspondence, characteristic varieties are identified with microsupports. That is, the statement ``$\cM$ quantizes $\mathrm{CC}(\cM)$" is read as ``$\mathrm{Sol}(\cM)$ quantizes $\SS(\mathrm{Sol}(\cM))$". The definition of SQv0 generalizes this to not only the constructible sheaves obtained as the image of the solution functor. In the real world, the Riemann--Hilbert correspondence does not exist, but sheaf quantization makes sense and is to some extent a topological concept. In this sense, sheaf quantization is more flexible.

\subsection{Weinstein case}
It was expected by several people that Fukaya categories of exact symplectic (especially, Weinstein) manifolds satisfy the locality (sheaf/cosheaf) property, and is called Kontsevich's conjecture (\cite{KontsevichSymp}). Nadler's program~\cite{NadlerCategoricalMorse} approaches this prediction by identifying the Fukaya category with a sheaf-theoretically constructed category. Although the original program relied on the arborealization conjecture~\cite{NadlerArboreal}\footnote{The conjecture itself has recent developments~\cite{AEN}.}, Ganatra--Pardon--Shende~\cite{MicrolocalMorse} has provided a solution that does not rely on the conjecture.

The sheaf-theoretic category was constructed in Shende~\cite{ShendeH-principle} and Nadler--Shende~\cite{nadlershende}.\footnote{It was expected that one can also construct the category by gluing, which has never been carried out in general to the author's knowledge.} We would like to outline their construction: Let $W$ be an exact symplectic manifold and $\lambda$ is a Liouville form. Let $W\times \bR$ be the contactification and embed it into $S^*\bR^N$ for sufficiently large $N$ as a contact manifold. Then the conormal bundle of the embedding is a symplectic bundle. Choose
one Lagrangian distribution $\tau$ in the normal bundle. This choice is called a polarization and is related to a choice of a stable Maslov structure. We denote the total space of $\tau$ by $\cU_\tau$. 

On the other hand, one can define the Kashiwara--Schapira stack over $T^*\bR^N$ by sheafifying the assignment $U\mapsto \Sh(\bR^N)/\Sh_{T^*\bR^N\bs U}(\bR^N)$. The category $\mu\mathrm{Sh}(\cU_\tau)$ defined as the global section of the stack over $\cU_\tau$.
It is proved using the $h$-principle that this does not depend on the choice of embedding, but only on the choice of the (stable) polarization.

From now on, suppose that $W$ is moreover Weinstein (see \S \ref{wrapped}). In the contact boundary $\partial W$, let us choose a Legendrian stop $S$ in it. Then this gives a subset $\bL\subset W$, called the relative skeleton.
The category we seek is the subcategory of $\mu\mathrm{Sh}(\cU_\tau)$ spanned by the objects whose (micro)support is in $\tau|_{\bL}$, where $\bL$ is considered as a subset of $W\times \bR$ through the embedding $W=W\times \{0\}\hookrightarrow W\times \bR$. We denote the resulting category by $\frakS\frakh(W, S)$.

Now we can state SQv0 for Weinstein manifolds.
\begin{thm}[Nadler--Shende~\cite{nadlershende}]
Let $L$ be a Lagrangian brane of $W$ ending at $S$. Then there is a fully faithful functor from the category of local systems $\mathrm{Loc}(L)$ to $\frakS\frakh(W, S)$.
\end{thm}

\begin{exa}[$A_n$-example, \cite{NadlerArboreal}]
Let us consider the case when $W=\bR^2_{x,y}$ equipped with $\lambda=xdy-ydx$. Then the skeleton is a point. One can embed $W\times \bR=\bR^3$ into $S^*\bR^2$. By the definition of $\frakS\frakh(W, \varnothing)$, this is a quotient category of the category of sheaves on $\bR^2$ whose microsupport is contained in a point at the contact boundary. Since the latter category is trivial, the former $\frakS\frakh(W, \varnothing)$ is so too.

Let us consider general two-dimensional cases. Let $W$ be a Weinstein surface. Then the skeleton locally looks like an $n$-valent tree. By an embedding of $W\times \bR$ into a contact boundary of $\bR^N$ for some $N$ and thickening by a choice of polarization, the $n$-valent tree becomes a Legendrian $n$-valent tree times a stabilization (i.e., times $\bR^k$ for some $k$). Along the stabilizing direction, sheaf theory does not change, so we will concentrate on the Legendrian $n$-valent tree.

Let us set
\begin{equation}
    f_i(x):=
    \begin{cases}
    &e^{-i/x} \text{ for $x\geq 0$}\\
    &0 \text{ for $x\leq 0$}
    \end{cases}
\end{equation}
for $i=0,..., n-1$. Let us consider the graph of $f_i$'s in $\bR^2_{x,y}$ and take the intersection with the unit ball $B:=\lc (x,y)\relmid x^2+y^2< 1\rc$. We denote it by $G_i$. This is a smooth submanifold in $B$. We take the $y$-positive part of the connected component of the conormal of $G_i$ and denote it by $N_i$. The union $\Lambda_n:=\bigcup_{i=0}^{n-1} N_i$ is an $n$-valent tree in the contact boundary $S^*B$. 

Consider the category
\begin{equation}
    \mu\mathrm{Sh}(S^*B):=\Sh(B)/\lc \cE \relmid \SS(\cE)\subset T^*_BB \rc.
    \end{equation}
Inside this category, we consider the full subcategory $\mu\mathrm{Sh}_{\Lambda_n}(S^*B)$ spanned by objects whose $\SS$ is contained in $\Lambda_n$ outside the zero-section. Then $\mu\mathrm{Sh}_{\Lambda_n}(S^*B)$ has a full exceptional collection $\la \cE_1,..., \cE_n \ra$ given by 
\begin{equation}
    \cE_i:=\bK_{\lc (x,y)\in \bR^2_{x,y}\relmid  y\geq f_i(x)\rc}.
\end{equation}
One can easily see that this exceptional collection gives an identification of $\mu\mathrm{Sh}_{\Lambda_n}(S^*B)$ with the representation category of $A_n$-quiver.

The whole category $\frakS\frakh(W, S)$ can be obtained as a gluing of these local categories using the sheaf property appeared in Kontsevich's conjecture.
\end{exa}

\begin{rem}
One can define the microlocal category for more general Maslov data by using some universal construction. See \cite{nadlershende} for details. 
\end{rem}

\subsection{Application of sheaf quantization version 0}
Here, of course, we do not mean applications of constructible sheaves to mathematics (which are too vast), but rather applications to symplectic topology and related areas.

The first application to symplectic topology is to Arnold's nearby Lagrangian conjecture \cite{NadlerMicrolocalbranes}. The nearby Lagrangian conjecture states that any compact exact Lagrangian submanifold in $T^*M$ is Hamiltonian isotopic to the zero-section. By the invariance of Lagrangian intersection Floer homology over a field under Hamiltonian isotopies, the conjecture implies that the endomorphism ring of a compact exact Lagrangian in $T^*M$ is isomorphic to the homology of $M$. Nadler proved this isomorphism by observing any compact exact Lagrangian corresponds to the rank 1 constant sheaf $\bK_M$ under the Nadler--Zaslow equivalence (a version of Theorem~\ref{theorem:NZGPS}).\footnote{Almost at the same time, Fukaya--Seidel--Smith proved a similar result without using sheaf theory~\cite{FSS}}.

The Nadler--Zaslow equivalence motivates an approach to homological mirror symmetry via sheaves. One of the main obstacles to prove HMS is the difficulty to compute Fukaya category. If we can interpret a Fukaya category as a sheaf category, the latter is usually easier to compute. Bondal~\cite{Bondal} and Fang--Liu--Treumann--Zaslow~\cite{FLTZ1, FLTZ2} initiated such an approach for the case of mirror symmetry for toric varieties. They define certain Lagrangians in the cotangent bundles of tori and Euclidean spaces and conjectured that the category of sheaves microsupported in the Lagrangians is equivalent to derived category of (usual/equivariant) coherent sheaves over toric varieties. They proved the equivariant case, and the non-equivariant case (extended to include the non-proper case) was later proved in \cite{kuwCCC} and the interpretation of the Lagrangian through mirror LG-models was given by \cite{GammageShende, Zhou}. 

\begin{exa}
The category of Example~\ref{exampleofsheafLagcorres2} is equivalent to the derived category of coherent sheaves over $\bP^1$. In fact, the objects $\bK_0, \bK_{(0,1)}[1]$ forms a full exceptional collection and can be identified with the collection $\la\cO_{\bP^1}, \cO_{\bP^1}(1)\ra$ of Beilinson collection~\cite{Beilinson}. The $\partial \bL$ can be obtained as a regular fiber of the mirror LG-model $W=z+\frac{1}{z}$.

Similarly, the category of Example~\ref{exampleofsheafLagcorres1} is equivalent to the derived category of coherent sheaves over the weighted projective stack $\bP^1(2,2)$.
\end{exa}

Other examples of HMS proved via sheaves include \cite{Nadlerwrapped, NadlerLG, GammageShende}... Now the complete list is huge and cannot be included here. Furthermore, as mentioned in the introduction, there are many other interesting applications of the sheaf-theoretic interpretation of Fukaya category.

\subsection{Coisotropic branes}
Fukaya category is a mathematical realization of the category of boundary conditions of A-model in physics (A-brane category). However, in Kapustin--Orlov~\cite{KapustinOrlov}, it was argued that the category of A-branes should contain such a brane that has support on a coisotropic submanifold. In other words, Fukaya category is only a subcategory of the A-brane category.\footnote{It might be possible that the subcategory is category-equivalent to the whole.} Furthermore, Gukov--Witten~\cite{gukovwitten} and Kapustin--Witten~\cite{KapustinWitten} pointed out that it is useful in mathematics (representation theory and the Langlands program) to consider coisotropic branes. However, Lagrangian intersection Floer theory for coisotropic submanifolds has not been developed yet. Using microlocal sheaf theory instead appears to be a good idea.

\begin{dfn}
Let $(X, \omega)$ be a symplectic manifold. Let $\frakC(X, \omega)$ be the set of coisotropic submanifolds of $X$. Then $\frakC$ is a coisotropic category of $(X, \omega)$ if the map $\frakC$ from the set $\mathrm{Ob}(\cC)$ to $\frakC(X, \omega)$ of the objects of $\cC$ is specified.
\end{dfn}

Let $M$ be a real manifold and $\bK$ a field. Let $\Sh(M)$ be the derived category of sheaves on $M$. By the Kashiwara--Schapira theorem~\cite{KashiwaraSchapira}, microsupport $\SS$ of sheaves are coisotropic. That is, the map $\SS$ equips $\Sh(M)$ with a coisotropic category structure. The category $\Sh(M)$ contains $\Sh^c(M)$ and $\Sh^w_{\Lambda}(M)$, that is, it contains Fukaya categories by versions of Theorem~\ref{theorem:NZGPS}. Therefore, it is natural to consider general objects of $\Sh(M)$ as coisotropic A-branes.

\begin{exa}[Canonical coisotropic brane]Let us recall the geometric Langlands conjecture by Beilinson--Drinfeld~\cite{BeilinsonDrinfeld}. Let us fix a Riemann surface $C$ and a reductive group $G$. We denote the derived category of quasi-coherent sheaves on the moduli $\cM(C, G)$ of $G$-local systems on $C$ by $D\mathrm{Qcoh}(\cM(C, G))$. We denote the derived category of $\cD$-modules on the moduli $\mathrm{Bun}_{{}^LG}(C)$ of principal $^{L}G$-bundles by $D(\cD_{\mathrm{Bun}_{{}^LG}(C))})$. Here ${}^LG$ is the Langlands dual of $G$. Then the geometric Langlands conjecture states that these categories are canonically equivalent.\footnote{This conjectures is known as a ``best-hope" version. A modified conjecture is proposed in \cite{ArinkinGaitsgory}.}

Kapustin--Witten~\cite{KapustinWitten} proposed to view this conjecture as a homological mirror symmetry between $\cM(C, G)$ and $T^*\mathrm{Bun}_{{}^LG}(C)$. Then it is necessary to relate the category of $\cD$-modules on $\mathrm{Bun}_{{}^LG}(C)$ to the A-brane category of $T^*\mathrm{Bun}_{{}^LG}(C)$. 

From now on let us set $M:=\mathrm{Bun}_{{}^LG}(C)$, to indicate the following arguments are valid for complex manifolds.
For regular holonomic $\cD$-modules, such a relation can be obtained by combining the Nadler--Zaslow equivalence and the Riemann--Hilbert correspondence \footnote{Rather, understanding Kapustin--Witten may have been one of Nadler--Zaslow's motivations.}:
\begin{equation}
    D^b_{rh}(\cD_M)\hookrightarrow \Sh^c(M)\cong \mathrm{Fuk}(T^*M).
\end{equation}
Since Beilinson--Drinfeld's geometric Langlands conjecture deals with general $\cD$-modules, we should avoid the restriction to regular holonomic ones.

Kapustin--Witten's argument proceeds as follows. First, they find an A-brane called the canonical coisotropic brane supported on the whole space $T^*M$. This brane has the property that its endomorphism ring is isomorphic to the ring of differential operators of the base space $M$. By the Yoneda pairing, we obtain a functor from the category of the A-branes of $T^*M$ to the category of the $\cD_M$-modules.

The sheaf corresponding to the canonical coisotropic brane is not in $\Sh^c(M)$ nor $\Sh^{wc}_{\bL}(M)$. However, there is an object in $\Sh(M)$ that seems like a canonical coisotropic brane.

First, consider $M$ as a real manifold. Let $C^\infty_M$ denote the sheaf of $C^\infty$-functions on $M$. Peetre's theorem~\cite{peetre} asserts that $\cEnd(C^\infty_M)$ is isomorphic to the sheaf of rings of $C^\infty$-differential on $M$. Also, $\SS(C^\infty_M)=T^*M$. Hence $C^\infty_M$ approximates the canonical coisotropic brane to some extent, but the endomorphism ring is not the ring of holomorphic/algebraic differential operators.

Peetre's theorem in the holomorphic category was conjectured by Sato and proved by Ishimura \cite{Ishimura} and explained in a sophisticated term in \cite{ProsmansSchneiders}. The claim is as follows. Let $M$ be a complex manifold. Let $\cO_M$ be the sheaf of holomorphic functions. Instead of $\Sh(M)$, consider the category of Ind-Banach space-valued sheaves $\Sh_{\mathrm{IB}}(M)$. One can view $\cO_M$ as an object of $\Sh_{\mathrm{IB}}(M)$. In this category $\cEnd(\cO_M)$ is isomorphic to $\cD^\infty_M$, the ring of differential operators with possibly infinite order. Despite the difference between $\cD_M$ and $\cD_M^\infty$, they are very close, for example, there is no difference between regular holonomic $\cD^\infty$-modules and regular holonomic $\cD$-modules where some important $\cD$-modules live in this class (e.g., \cite{EFK}). Hence, we can consider $\cO_M$ as a good approximation of the canonical coisotropic brane.\footnote{To have a better approximation, one could try to prove Ishimura-type theorems in the setting of ind-sheaves, where some tame versions of structure sheaf are available~\cite{KSind}.}

There is also a relevance to the regular Riemann--Hilbert correspondence. If we consider $\cO_M$ as the canonical coisotropic brane, the Yoneda module appeared in Kapustin--Witten's proposal restricts to
\begin{equation}
\cHom(-, \cO_M)\colon \Sh^c(M)\rightarrow D^b(\cD^\infty_M),
\end{equation}
which is the converse to the solution functor~\cite{KashiwaraRH, Mebkhout, ProsmansSchneiders}.

It is also consistent with the recent proposal of Gaiotto--Witten~\cite{gaiotto2021probing} as follows: Let us consider the constant sheaf $\bK_M$, which is a brane supported on the zero-section of $T^*M$. Then $\Hom(\bK_M, \cO_M)=\Gamma(M, \cO_M)$ is the space of holomorphic functions. This is the state space of geometric quantization when the fiber polarization is taken.
\end{exa}

\begin{rem}
Although tautological, if one takes the category of the $\cD$-modules as a mathematical definition of the A-brane category (in the Langlands setup), then one can consider $\cD$ as the canonical coisotropic brane (e.g. \cite{Frenkelconformalblock}). Then $\Hom(\cD, \cD)\cong \cD$, and the characteristic variety of $\cD$ is $T^*M$. On the other hand, an advantage of using sheaves (i.e., instead of $\cD$-modules) is that the object is topological. For example, the GKS-functor can act on it.
\end{rem}

\begin{rem}
Coisotropic subsets can also be defined for Poisson manifolds. Thus, the notion of coisotropic category can also be defined for Poisson manifolds. Attempts to consider A-branes on Poisson manifolds were made by Gualtieri~\cite{Gualtieri} and Bischoff--Gualtieri~\cite{BischoffGualtieri}. It is an interesting challenge to consider the problem via microlocal sheaf theory.
\end{rem}

\section{Sheaf quantization version 1}
\subsection{Tamarkin's trick}
The problem of SQv0 is that the image of the structural morphism of the Lagrangian category is contained in the set of conic subsets in $\cL(T^*M, \omega)$, since microsupports are always conic. In other words, Lagrangian submanifolds that can be treated from a sheaf-theoretic point of view are limited to only conic ones.

One can criticize the above claim as follows: Since Nadler--Zaslow-type theorems hold and a Fukaya category does include non-conic Lagrangians, there are no problems like the above. OK, we accept this objection.  
But, this approach also has some problems. 
\begin{enumerate}
    \item To treat a given Lagrangian, we have to first make it conic. In this process, geometrically concrete information (e.g., the number of intersection points) will be lost.
\item In other words, microlocal sheaf theory can only treat symplectic geometric data that are invariant under the isomorphisms in Fukaya categories, in particular, Hamiltonian diffeomorphisms. 
\end{enumerate}

Tamarkin's trick partially solves these problems. Let $M$ be a real manifold. The standard Liouville form is of the form $\sum_{i=1}^n\xi_idx_i$ in local coordinates, where $\{x_i\}$ is a local coordinate of $M$ and $\xi_i$ are the corresponding cotangent coordinates.

A Lagrangian submanifold $L$ of $T^*M$ is said to be exact if $\lambda|_L$ is exact.\footnote{Note that $\lambda|_L$ being closed is equivalent to the definition of Lagrangian.} Tamarkin's trick proceeds as follows: First, take a primitive $f\colon L\rightarrow \bR$ of $\lambda|_L$.

Let $\bR_t$ be the real line with standard coordinate $t$. Let $T^*\bR_t$ have the trivialization $T^*\bR_t\cong \bR_t\times \bR$ by the global section $dt$. The corresponding cotangent coordinate is written as $\tau$. Write $T^*_{\tau>0}\bR_t$ for the open set defined by $\tau>0$ in $T^*\bR_t$. Then
\begin{equation}
    L_f:=\lc (x, \xi, t, \tau)\relmid (x, \xi)\in T^*_xM, \tau>0, (x, \xi/\tau)\in L, t=-f(x, \xi/\tau)\rc
\end{equation}
is a Lagrangian submanifold of $T^*M\times T^*\bR_t$. We can recover $L$ from $f$ by the twisted projection
\begin{equation}
    \rho\colon T^*M\times T^*_{\tau>0}\bR_t\rightarrow T^*M; (x, \xi, t, \tau)\rightarrow (x, \xi/\tau)
\end{equation}
as $\rho(L_f)=L$.

The point here is that $L$ is not conical, so there are no hope to relate microsupport of sheaves to $L$, but $L_f$ is conical. The basic idea is to consider a sheaf microsupported on $L_f$. However, from the definition, $L_f$ is contained in $\tau>0$. The definition of microsupport tells us that $\SS(\cE)\cap T^*_{M\times \bR}M\times \bR=\varnothing$ implies $\cE=0$. In order to have a meaningful microsupport only in $\tau>0$, we localize the category:
\begin{equation}
    \Sh_{\tau >0}(M\times \bR):=\Sh(M\times \bR)/\lc \cE\in \Sh(M\times \bR)\relmid \SS(\cE)\subset \lc \tau\leq 0\rc\rc.
\end{equation}
This category is sometimes called the Tamarkin category. For $\cE\in \Sh_{\tau>0}(M\times \bR)$, $\SS(\cE)\cap \{\tau>0\}$ is an invariant concept under isomorphisms in $\Sh_{\tau>0}(M\times \bR)$.

\begin{rem}
There exists a canonical embedding of $\Sh_{\tau >0}(M\times \bR)$ into $\Sh(M\times \bR)$. In this sense, $\Sh_{\tau >0}(M\times \bR)$ can be considered as a category of sheaves.
\end{rem}

\begin{dfn}
We say an object $\cE\in \Sh_{\tau>0}(M\times \bR)$ is a sheaf quantization version 1 (SQv1 in short) of an exact Lagrangian $L$ if  $\SS(\cE)\cap \{\tau>0\}=L_f$.

We sometimes moreover assume that the microstalk condition (or the type condition in the terminology of \cite{KashiwaraSchapira}). See Subsection 4.4 for more details. 

We say an object $\cE$ is an SQv1 if it is in the subcategory of $\Sh_{\tau>0}(M\times \bR)$ generated by SQv1 of exact Lagrangians.
\end{dfn}

It is reasonable to define the non-conic version of microsupport using $\rho$:
\begin{dfn}[Non-conic microsupport]
Let $\cE$ be an object of $\Sh_{\tau>0}(M\times \bR)$. We set
\begin{equation}
    \mu\mathrm{supp}(\cE):=\rho(\SS(\cE)\cap \{\tau>0\}).
\end{equation}
\end{dfn}
Then a sheaf quantization $\cE$ of $L$ satisfies $\mu\mathrm{supp}(\cE)=L$. The category of SQv1 has a structure of a Lagrangian category by $\mu\mathrm{supp}$.

The following existence theorem for SQv1 was first stated by Viterbo~\cite{viterbo} for compact Lagrangians, proved by Guillermou~\cite{Guillermou} for compact exact Lagrangians and Jin--Treumann~\cite{JT} for non-compact exact Lagrangians.
\begin{thm}[Guillermou, Jin--Treumann] Any ``good" exact Lagrangian submanifold with a brane structure admits a sheaf quantization.
\end{thm}

\begin{rem}
Viterbo stated the above theorem in compact case in \cite{viterbo} and indicated to use Floer theory to prove it. The proof along this line is later provided in \cite{Viterbo2}. This Floer-theoretic proof is important, since one can directly relate Floer-theoretic consequences with sheaf-theoretic ones. 

Another important point of view provided by his proof (and also in Tamarkin's original construction~\cite{Tam}) is that sheaf quantization is, in some sense, a generalization of the theory of generating functions~\cite{Viterbo3}.
\end{rem}

The ``goodness" of Lagrangian submanifold asserts a certain good behavior of $L$ at $\infty$ of $T^*M$ (see \cite{JT} for detail, where goodness=tameness$+$lower exactness). Brane structures will be explained in Section~\ref{section:brane}. Guillermou moreover sheaf-theoretically proved that any compact exact Lagrangian submanifold in $T^*M$ admits a brane structure. Hence any compact exact Lagrangian submanifold in $T^*M$ admits a sheaf quantization.

\begin{exa}\label{example:Sqv1}
\begin{enumerate}
    \item Let $L$ be a conic Lagrangian. Let $\cE$ be an SQv0 of $L$. We set $\bR_{t\geq 0}:=\lc t\in \bR\relmid t\geq 0 \rc$. The inclusion $\iota\colon \bR_{t\geq 0}\hookrightarrow \bR$. By abuse of notation, we denote the push-forward of $\bK_{t\geq 0}$ along $\iota$ by $\bK_{t\geq 0}$. The exterior tensor $\cE\boxtimes \bK_{t\geq 0}\in \Sh_{\tau>0}(M\times \bR)$ defines an SQv1 of $L$.
    \item Let $f\colon M\rightarrow \bR$ be a $C^1$-function. Let $L:=\mathrm{Graph}(df)$ be a Lagrangian submanifold of $T^*M$. We denote the push-forward of the rank 1 constant sheaf on $\lc(x,t)\in M\times \bR_t\relmid t\geq -f(x) \rc$ along the inclusion to $M\times \bR_t$ by $\bK_{t\geq -f(x)}$. Then this is a sheaf quantization of $L$. See also Figure~\ref{Figure:sqv1}.
    \begin{figure}
  \begin{minipage}[b]{0.45\linewidth}
    \centering
\begin{tikzpicture}
    \begin{scope}[gray]
    \draw[->](-2,0) -- (2,0);
    \draw[->](0,-1) -- (0,2.1);
  \end{scope}
\draw[-,thick] (-2,2)-- (1,-1);
\draw (0.5,-0.5) node[above] {$L$};
\draw (1,1.8) node[right] {$T^*\bR$};
\draw[-] (1,1.5)-- (1,2);
\draw[-] (1.5,1.5)-- (1,1.5);
\end{tikzpicture}
\subcaption{$L=\{\xi=d(-x^2/2)\}$}
  \end{minipage}
  \begin{minipage}[b]{0.45\linewidth}
      \centering
  \begin{tikzpicture}
  \draw (0,2.1) node[above] {$\bR_t$};
  \draw (2,0) node[right] {$\bR$};
  \filldraw[fill=lightgray, thick,domain=-1.4:1.4] plot(\x,\x*\x);
      \begin{scope}[gray]
    \draw[->](-2,0) -- (2,0);
    \draw[->](0,-1) -- (0,2.1);
  \end{scope}
  \end{tikzpicture}
\subcaption{An SQv1 of $L$}
  \end{minipage}
  \caption{Example~\ref{example:Sqv1}}\label{Figure:sqv1}
\end{figure}
\end{enumerate}
\end{exa}

Stronger results on sheaf quantizations of graphs of Hamiltonian isomorphisms were obtained earlier than the above theorem. 
We say a function $f$ defined on $T^*M\bs T^*_MM$ is conic if it satisfies $f(x, a\xi)=a f(x, \xi)$ for $a>0$. 

If $H$ is a compactly-supported function, we can consider the conification $\widehat{H}$ on $T^*M\times T^*_{\tau\neq 0}\bR$ by 
\begin{equation}
    \widehat{H}:=\tau H(x, \xi/\tau).
\end{equation}
In the below theorem, the second statement is based on the first result and this conification procedure.
\begin{thm}[Guillermou--Kashiwara--Schapira~\cite{GKS}]
\begin{enumerate}
    \item Let $\widehat{H}$ be a conic Hamiltonian function (not necessarily obtained as a conification) on $T^*M\bs T^*_MM$. Let $\phi_{\widehat{H}}$ be the associated conic time-1 Hamiltonian isomorphism. Then there exists a category equivalence $\Phi_{\widehat{H}}\colon \Sh(M)\rightarrow \Sh(M)$ satisfying $\phi_H(\SS(\cE))\bs T^*_MM)=\SS(\Phi(\cE))\bs T^*_MM$.
    \item Let $H$ be a compactly-supported Hamiltonian function on $T^*M$. Let $\phi_H$ be the associated time-1 Hamiltonian isomorphism. Then there exists a category equivalence $\Phi_H\colon \Sh_{\tau>0}(M\times \bR)\rightarrow \Sh_{\tau>0}(M\times \bR)$ satisfying $\phi_H(\mu\mathrm{supp}(\cE))=\mu\mathrm{supp}(\Phi(\cE))$. The equivalence is given by an integral transform with a kernel sheaf-quantizing the graph of $\Phi$.
\end{enumerate}
\end{thm}
The obtained functor is sometimes called GKS-functor, GKS-equivalence, etc. An earlier primitive version of the above theorem can be found in Tamarkin~\cite{Tam}

\begin{exa}\label{ex:GKS}
We consider the simplest case when $M=\bR$. Let us denote the coordinate function on $\bR$ by $x$ and the cotangent coordinate by $\xi$. Consider the Hamiltonian function $|\xi|$. Then the Hamiltonian vector field is given by $-\partial_x$ on $T^{*,+}\bR$ and $\partial_x$ on $T^{*,-}\bR$. The GKS-functor at time $t$ is given by the convolution with $\bK_{\lc(x,y)\relmid |x-y|\leq t\rc}$. In particular, the kernel at $t=0$ is the constant sheaf on the diagonal as expected. For $t<0$, the kernel is given by $\bK_{\lc(x,y)\relmid |x-y|<-t\rc}[1]$.
\begin{figure}
\begin{minipage}[b]{0.45\linewidth}
    \centering
      \begin{tikzpicture}
    \begin{scope}[gray]
    \draw[->](-2,0) -- (2,0);
    \draw[->](0,-2) -- (0,2);
  \end{scope}
  \draw (0,2.1) node[above] {$x_2$};
  \draw (2,0) node[right] {$x_1$};
  \draw[thick] (-2,-2) -- (2,2);
 \end{tikzpicture}
 \subcaption{$t=0$}
\end{minipage}
\begin{minipage}[b]{0.45\linewidth}
\centering
      \begin{tikzpicture}
  \draw (0,2.1) node[above] {$x_2$};
  \fill[fill=lightgray] (-1,-2) -- (2,1) -- (2,2) -- (1,2) -- (-2, -1) -- (-2, -2) -- (-1,-2);
  \draw (2,0) node[right] {$x_1$};
  \draw[thick] (-1,-2) -- (2,1);
  \draw[thick] (-2,-1) -- (1,2);
      \begin{scope}[gray]
    \draw[->](-2,0) -- (2,0);
    \draw[->](0,-2) -- (0,2);
  \end{scope}
 \end{tikzpicture}
 \subcaption{$t=1$}
\end{minipage}
\caption{Example~\ref{ex:GKS}}
\end{figure}

\end{exa}

\begin{rem}
For the graph of a Hamiltonian isotopy, the brane data (\S \ref{section:brane}) is canonically determined from the canonical brane data of the diagonal, since the notion of brane data is homotopical. 
\end{rem}

\subsection{Shift and Energy}
To illustrate applications of SQv1, a few more words about the properties of the Tamarkin category. First, let $m\colon \bR_t\times \bR_t\rightarrow \bR_t$ be the addition morphism. The induced morphism $(M\times \bR_t)^2\rightarrow M\times M\times \bR_t$ will be written as $m$ as well. For $\cF, \cG\in \Sh(M\times \bR_t)$, we set $\cF\star \cG:=m_! (\cF\boxtimes \cG)$. This induces a product structure in $\Sh_{\tau>0}(M\times \bR_t)$. The unit object is $\bK_M\boxtimes \bK_{t\geq 0}$. For $c\in \bR$, the product by $\bK_M\boxtimes \bK_{t\geq c}$ is a translation by $c$. We will write this operation as $T_c$.

For $c\leq c'$, there is a canonical morphism $T^{c'-c}\colon \bK_{t\geq c}\rightarrow \bK_{t\geq c'}$ corresponding to $\id_{\bK_{t\geq c'}}$ under the adjunction $\Hom(\bK_{t\geq c}, \bK_{t\geq c'})\cong \Hom(\bK_{t\geq c'}, \bK_{t\geq c'})$. This satisfies $T^aT^b=T^{a+b}$ and induces $T^{c'-c}\colon T_c\cF\rightarrow T_{c'}\cF$ through the product. See Figure~\ref{fig:translation}.

\begin{figure}
\centering
  \begin{tikzpicture}
    \begin{scope}[gray]
    \draw[->](-2.1,-0.5) -- (2,-0.5);
    \draw[->](-2,-1) -- (-2,2.1);
  \end{scope}
  \draw (-2,2.1) node[above] {$\bR_t$};
  \draw (2,0) node[right] {$\bR$};
  \filldraw[fill=lightgray, thick,domain=-1.42:1.42] plot(\x,\x*\x);
  \filldraw[fill=gray, thick,domain=-1:1] plot(\x,\x*\x+1);
    \draw (0,0.6) node[below] {$\cE$};
 \draw (0,1.2) node[above] {$T_c\cE$};
\draw[->](0,0.6) -- (0, 1.2);
\draw (0,0.85) node[right] {$T^c$};
  \end{tikzpicture}
    \caption{$T^c$}
    \label{fig:translation}
\end{figure}

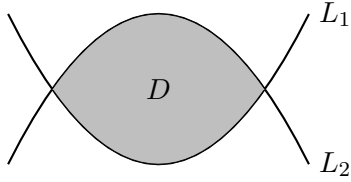
\begin{figure}
    \centering
      \begin{tikzpicture}
\draw[thick,domain=-2:2] plot(\x,0.5*\x*\x-1);
\draw[thick,domain=-2:2] plot(\x,-0.5*\x*\x+1);
\filldraw[fill=lightgray, domain=-1.42:1.42] plot(\x,0.5*\x*\x-1);
\filldraw[fill=lightgray, domain=-1.42:1.42] plot(\x,-0.5*\x*\x+1);
\draw (0,0) node {$D$};
\draw (2,1) node[right] {$L_1$};
\draw (2,-1) node[right] {$L_2$};
  \end{tikzpicture}
\caption{Lagrangian intersections}\label{fig:Lagintersection}
\end{figure}

Later, this will be interpreted as a Novikov ring action. To present an evidence, we would like to recall the definition of the differential of Hom-spaces of Fukaya category. For an intersection point $p\in \Hom(L_1, L_2)$, the differential of $p$ is given by the sum $\sum_{q\in L_1\cap L_2}\#\cM(L_1, L_2, p,q)q$, where the number is again counting certain holomorphic disks of the form Figure~\ref{fig:Lagintersection}. 

Let us consider the situation where there are exactly two points of intersections $p,q\in L_1\cap L_2$ and $dp=q$ holds, in particular, only one holomorphic disk exists. Then the hom-complex $\Hom(L_1, L_2)$ is acyclic. Nothing remains.

But, there exists a refinement. Instead of considering $dp=q$, we consider $dp=T^{\int_D\omega} q$, where $T$ is an indeterminate. As we shall introduce \S5, this is the Novikov parameter. In this refinement, $q$ is cohomologically 0 only after multiplied by $T^{\int_D\omega}$.

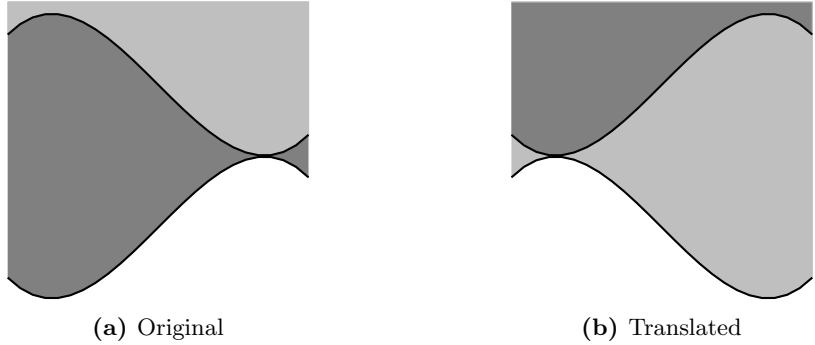
\begin{figure}
\begin{minipage}[b]{0.45\linewidth}
    \centering
         \begin{tikzpicture}
\filldraw[fill=gray,draw=lightgray,domain=-2:2] plot(\x,-0.5*\x*\x*\x/3+\x) -- (2,3) -- (-2,3);
\filldraw[fill=lightgray,draw=lightgray,domain=-2:2] plot(\x,0.5*\x*\x*\x/3-\x+1.9) -- (2,3) -- (-2,3);
\draw[thick,domain=-2:2] plot(\x,0.5*\x*\x*\x/3-\x+1.9);
\draw[thick,domain=-2:2] plot(\x,-0.5*\x*\x*\x/3+\x);
  \end{tikzpicture}
 \subcaption{Original}\label{fig:S12}
\end{minipage}
\begin{minipage}[b]{0.45\linewidth}
\centering
\begin{tikzpicture}
\filldraw[fill=lightgray,draw=lightgray,domain=-2:2] plot(\x,0.5*\x*\x*\x/3-\x) -- (2,3) -- (-2,3);
\filldraw[fill=gray,draw=lightgray,domain=-2:2] plot(\x,-0.5*\x*\x*\x/3+\x+1.9) -- (2,3) -- (-2,3);
\draw[thick,domain=-2:2] plot(\x,0.5*\x*\x*\x/3-\x);
\draw[thick,domain=-2:2] plot(\x,-0.5*\x*\x*\x/3+\x+1.9);
  \end{tikzpicture}
 \subcaption{Translated}\label{fig:S12c}
\end{minipage}
\caption{SQv1 of $L_1$ and $L_2$}
\end{figure}

A similar thing can be observed on the sheaf side. Let us choose primitives of $L_1$ and $L_2$ in Figure~\ref{fig:Lagintersection} like in the curves of Figure~\ref{fig:S12}. Then an SQv1 $S_1$ of $L_1$ is given by the constant sheaf supported on the light-gray region an SQv1 $S_2$ of $L_2$ is given by the constant sheaf supported on the light-gray on the light-gray and the gray regions.
We have a morphism $S_1\rightarrow S_2$ (with a possible shift of degree). By the composition with $T^c$, this morphism induces $S_1\rightarrow T_cS_2$ for any $c\geq 0$. When the picture becomes Figure~\ref{fig:S12c}, the composition will vanish. This is exactly when $c=\int_D\omega$. 

\subsection{Applications of sheaf quantization version 1}
We first describe Tamarkin's original work~\cite{Tam}. For an explanation, see for example Guillermou--Schapira~\cite{GuiS}. Note that a Hamiltonian isotopy has energy and drives primitive functions. To obtain a Hamitonian-isotopy-invariant category, on has to quotient out the torsion objects:
\begin{dfn}
We say $\cF\in \Sh_{\tau>0}(M\times \bR)$ is a torsion if there exists some $c\in \bR_{\geq 0}$ such that $T^c\colon \cF\rightarrow T_c\cF$ is a zero map.
\end{dfn}
Write $\cT(M)$ for the category obtained from $\Sh_{\tau>0}(M\times \bR)$ by quotienting out the torsion object. The following theorem holds.
\begin{thm}
Let $\Phi$ be the GKS-functor for a compactly supported Hamiltonian. Then for any sheaf $\cE$, we have $\cE\cong \Phi(\cE)$ in $\cT(M)$.
\end{thm}

As a corollary of the above theorem, we obtain the following result, called Tamarkin's non-displaceablity theorem.
\begin{crl}
Suppose $\cF, \cG\in \cT(M)$ satisfies $\Hom(\cF, \cG)\neq 0$. Then $\mu\mathrm{supp}(\cF)$ cannot be moved to not intersect with $\mu\mathrm{supp}(\cG)$ by any compactly-supported Hamiltonian isotopy.
\end{crl}
\begin{proof}
Suppose there exists a compactly supported Hamiltonian such that the associated Hamiltonian isomorphism displaces $\mu\mathrm{supp}(\cG)$ from $\mu\mathrm{supp}(\cF)$. Let us denote the associated GKS-functor by $\Phi$. Then, by a microsupport estimate, one can conclude that $\Hom(\cF, \Phi(\cG))=0$. Using the above theorem, we can complete the proof.
\end{proof}

As we have already mentioned, Guillermou~\cite{Guillermou} showed that any compact exact Lagrangian in $T^*M$ can be sheaf-quantized and it can be used to conclude that any compactly-supported Hamiltonian isotopy cannot displace a compact exact Lagrangian from itself. He also used it to prove a version of the Arnold conjecture. Guillermou's sheaf-theoretic study of symplectic geometry including the above result is collected in \cite{Guillermoucollection}. Ike~\cite{IkeLagrangianintersection} compares Guillermou's sheaf quantization with Tamarkin's result in detail to show an estimate of the number of Lagrangian intersections sheaf-theoretically.

Although it is not a proper application, the idea of adding an extra variable $\bR_t$ has been found useful in the Riemann--Hilbert correspondence for possibly irregular $\cD$-modules by D'Agnolo--Kashiwara~\cite{DK}. Let us briefly explain this. For the regular case, the monodromy of solutions is enough to recover the original differential equation. However, for the irregular case, the monodromy is not enough and we have to record the Stokes data, which can be considered as ``sectorial monodromy around the singularities" (see, for example, \cite{Wasow}). When the base manifold is of $\dim_\bC=1$, the equivalence between the Stokes data and irregular differential equations was established by Deligne--Malgrange--Sibuya. One of the key ingredients is the formal classification result of irregular singularities due to Hukuhara--Levelt--Turrittin, which is missing in higher-dimensional cases for many years. Based on the recent establishment of HLT-theorem in higher dimensions by Sabbah--Kedlaya--Mochizuki~\cite{SabbahIrreg, Kedlaya, WildTwistor}, D'Agnolo--Kashiwara proved the RH-correspondence of holonomic $\cD$-modules. One of the key ideas in their work is that one can encode Stokes data in the extra variable $\bR_t$. How? The formal type of irregular singularity is given by a set of Laurent functions modulo regular functions. Then the extra $\bR_t$ is used as the codomain of the real part of the formal type. But this is not well-defined, since there is an ambiguity to add some regular (i.e., bounded) functions. To ignore this ambiguity, D'Agnolo--Kashiwara ignores the finite translation along $\bR$.
This is quite similar to the goal of Tamarkin who achieved it by quotienting out the torsion objects, but D'Agnolo--Kashiwara instead takes a way taking the limit of translations, namely, $\lim_{c\rightarrow \infty} T_c$. This limit does not make sense as it is, so one has to compactify $\bR_t$ and consider ind-limits, then one will arrive at the notion of enhanced ind-sheaves, which is the target category of their RH correspondence. For more details, one should refer to the original article~\cite{DK} or the survey~\cite{KSregularandirregular}.

This relation of Stokes data and the extra $\bR_t$-direction can be considered as a degenerate version of the use in \cite{kuwagaki2020sheaf, kuwagakihRH}. We will briefly mention this later.

\begin{rem} Somehow the direction $\bR_t$ is not geometric, we eventually have to ``squash" this direction. Tamarkin did it by quotienting out the torsion objects, Polesello--Schapira~\cite{PS} did it by considering $t$-independency of functions, D'Agnolo--Kashiwara did it by taking the limit of translations, the author did it by considering equivariant sheaves in \cite{kuwagaki2020sheaf}, and by considering the finite Novikov ring in \cite{KuwIrreg}.
\end{rem}

As already mentioned, it was pointed out by Tamarkin that the added parameter $\bR_t$ is related to symplectic energy (i.e., Novikov exponent): instead of quotienting out the torsion objects as Tamarkin did, it was hoped that by looking at the direction $\bR_t$ in detail, quantitative conclusions can be drawn. This was actually done by Asano--Ike~\cite{AsanoIkepersitence}. Their notion of interleaving distance is based on the work of Kashiwara--Schapira~\cite{KashiwaraSchapirapersistent} on a sheaf-theoretic treatment of persistent homology.

\begin{dfn}[Interleaving distance]
For $\cE, \cF\in \Sh_{\tau>0}(M\times \bR)$, we say $\cF$ is $(a,b)$-interleaved for $(a,b)\in \bR^2_{\geq 0}$ if there exist morphisms $\alpha, \alpha'\colon \cE\rightarrow T_a\cF$ and $\beta, \beta'\colon \cF\rightarrow T_b\cE$ such that $T_a(\beta)\circ \alpha=T^{a+b}$ (resp. $T_b(\alpha')\circ \beta'=T^{a+b}$) as morphisms from $\cE$ to $T_{a+b}(\cE)$ (resp. from $\cF$ to $T_{a+b}(\cF)$). We set
\begin{equation}
    d(\cE, \cF):=\inf\lc a+b\relmid \text{$\cE$ and $\cF$ are $(a,b)$-interleaved}\rc.
\end{equation}
This defines a pseudo-metric on the set of objects of $\Sh_{\tau>0}(M\times \bR)$.
\end{dfn}
It follows from the following theorem that this can be used to estimate symplectic energy. Let $H$ be a function on $X\times I$.
\begin{equation}
    ||H||_{\mathrm{osc}}:=\int_0^1\lb \max_{x\in X}H(x,s)-\min_{x\in X}H(x,s)\rb ds.
\end{equation}
\begin{thm}[Asano--Ike~\cite{AsanoIkepersitence}]Let $H$ be a compactly-supported Hamilton function, $\Phi$ be the associated GKS-functor, and $\cF\in \Sh_{\tau>0}(M\times \bR)$. In this situation, we have
\begin{equation}
    d(\cF, \Phi(\cF))\leq ||H||_{\mathrm{osc}}.
\end{equation}
\end{thm}
Asano--Ike used this to evaluate displacement energy. For further applications, see e.g., \cite{Li, Zhang}.

\begin{exa}\label{4.12}
Let us go back to the example of $M=S^1=\bR/\bZ$. Recall that we have ``SQv0" of $\SS(\bK_{(0,1)})=\SS(\bK_{[0,1]})$. We consider the corresponding SQv1 $\bK_{(0,1)}\boxtimes \bK_{t\geq 0}\in \Sh_{\tau>0}(M\times \bR)$ and $\bK_{[0,1]}\boxtimes \bK_{t\geq 0}\in \Sh_{\tau>0}(M\times \bR)$. For $a\in (0, 1/2)$, we would like to consider the following exact Lagrangian in $T^*S^1$ (see Figuer~\ref{fig:La}):
\begin{equation}
\begin{split}
        L_a=&p_{T^*S^1}([a, 1-a]\cup \lc (x, \xi)\relmid (x-a)^2+(\xi-a)^2=a^2, x-a\leq 0, \xi-a\leq 0\rc \\
        &\cup \lc (x, \xi)\relmid (x+a)^2+(\xi+a)^2=a^2, x+a\geq 0, \xi+a\geq 0\rc\cup \lc (0, \xi)\relmid |\xi|>a \rc)
\end{split}
\end{equation}
where $p_{T^*S^1}\colon \bR^2\rightarrow T^*S^1$ is the standard projection.

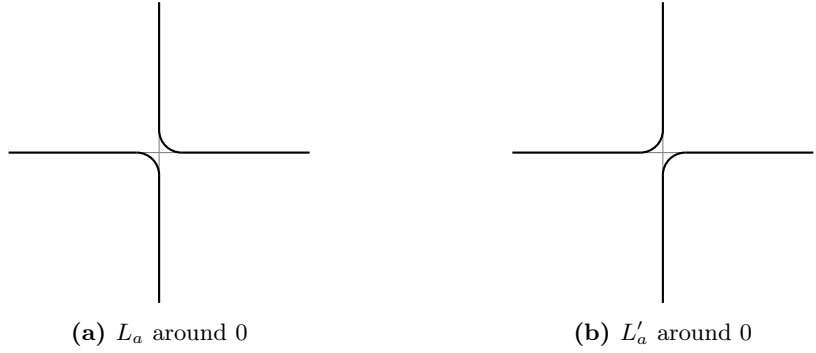
\begin{figure}
\begin{minipage}[b]{0.45\linewidth}
    \centering
  \begin{tikzpicture}
    \begin{scope}[gray]
    \draw[-](-2,0) -- (2,0);
    \draw[-](0,-2) -- (0,2);
  \end{scope}
\draw[-,thick](2,0) -- (0.3,0);
\draw[-,thick](0,2) -- (0,0.3);
\draw[-,thick](-2,0) -- (-0.3,0);
\draw[-,thick](0,-2) -- (0,-0.3);
\draw[thick](0,0.3) arc (180:270:0.3cm);
\draw[thick](-0.3,0) arc (90:0:0.3cm);
  \end{tikzpicture}
 \subcaption{$L_a$ around $0$}\label{fig:La}
\end{minipage}
\begin{minipage}[b]{0.45\linewidth}
\centering
  \begin{tikzpicture}
    \begin{scope}[gray]
    \draw[-](-2,0) -- (2,0);
    \draw[-](0,-2) -- (0,2);
  \end{scope}
\draw[-,thick](2,0) -- (0.3,0);
\draw[-,thick](0,2) -- (0,0.3);
\draw[-,thick](-2,0) -- (-0.3,0);
\draw[-,thick](0,-2) -- (0,-0.3);
\draw[thick](0,0.3) arc (0:-90:0.3cm);
\draw[thick](0.3,0) arc (90:180:0.3cm);
  \end{tikzpicture}
 \subcaption{$L_a'$ around $0$}\label{fig:La'}
\end{minipage}
\caption{Lagrangians of Example~\ref{4.12}}
\end{figure}

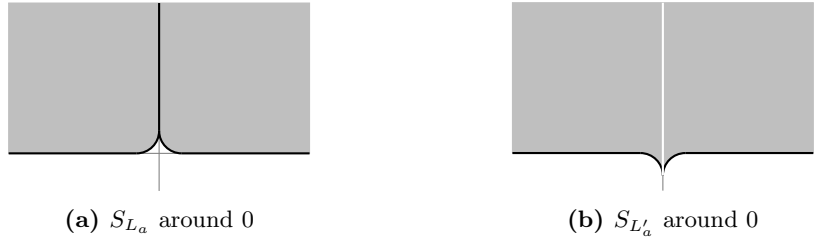
\begin{figure}
\begin{minipage}[b]{0.45\linewidth}
    \centering
    \begin{tikzpicture}
    \begin{scope}[gray]
    \draw[-](-2,0) -- (2,0);
    \draw[-](0,-0.5) -- (0,2);
  \end{scope}
\filldraw[fill=lightgray,draw=lightgray](0,0.3) arc (180:270:0.3cm) -- (2,0)--(2,2)--(0,2);
\filldraw[fill=lightgray,draw=lightgray](0,0.3) arc (0:-90:0.3cm)-- (-2,0)--(-2,2)--(0,2);
\draw[-,thick](2,0) -- (0.3,0);
\draw[-,thick](0,2) -- (0,0.3);
\draw[-,thick](-2,0) -- (-0.3,0);
\draw[thick](0,0.3) arc (0:-90:0.3cm);
\draw[thick](0,0.3) arc (180:270:0.3cm);
  \end{tikzpicture}
 \subcaption{$S_{L_a}$ around $0$}\label{fig:SLa}
\end{minipage}
\begin{minipage}[b]{0.45\linewidth}
\centering
  \begin{tikzpicture}
    \begin{scope}[gray]
    \draw[-](-2,0) -- (2,0);
    \draw[-](0,-0.5) -- (0,2);
  \end{scope}
\filldraw[fill=lightgray,draw=lightgray](0.3,0) arc (90:180:0.3cm) -- (0,2)--(2,2)--(2,0);
\filldraw[fill=lightgray,draw=lightgray](-0.3,0) arc (90:0:0.3cm)-- (0,2)--(-2,2)--(-2,0);
\draw[-,thick](2,0) -- (0.3,0);
\draw[-,thick](-2,0) -- (-0.3,0);
\draw[thick](0.3,0) arc (90:180:0.3cm);
\draw[thick](-0.3,0) arc (90:0:0.3cm);
\draw[-,thick,white](0,2) -- (0,-0.3);
  \end{tikzpicture}
 \subcaption{$S_{L_a'}$ around $0$}\label{fig:SLa'}
\end{minipage}
\caption{Sheaf quantizations of Example~\ref{4.12}}
\end{figure}

This Lagrangian admits a sheaf quantization (Figure~\ref{fig:SLa}):
\begin{equation}
\begin{split}
    S_{L_a}&:=p_!\bK_{D_a}\\
    D_a&:=\lc (x, t)\relmid t\geq 0, x\in (a, 1-a)\rc\cup \lc (x, t)\relmid x\in [0, a], t\geq \int_x^a\lb -\sqrt{a^2-(x-a)^2}+a\rb dx\rc\\
    &\cup \lc (x, t)\relmid x\in [1-a, 1], t\geq \int_{1-a}^x\lb -\sqrt{a^2-(x-1+a)^2}+a\rb dx\rc
\end{split}
\end{equation}
where $p$ is the projection $p\colon \bR\times \bR_t\rightarrow S^1\times \bR_t$.
Then one can see that $S_{L_a}\rightarrow \bK_{[0,1]}\boxtimes \bK_{t\geq 0}$ as $a\rightarrow 0$ with respect to the interleaving distance.

On the other hand, if we set 
\begin{equation}
\begin{split}
    L_a'=&p_{T^*S^1}([a, 1-a]\cup \lc (x, \xi)\relmid (x+a)^2+(\xi-a)^2=a^2, x+a\geq 0, \xi-a\leq 0\rc\\
        &\cup  \lc (x, \xi)\relmid (x-a)^2+(\xi+a)^2=a^2, x-a\leq 0, \xi+a\geq 0\rc\cup \lc (0, \xi)\relmid |\xi|>a \rc)
\end{split}
\end{equation}
(Figure~\ref{fig:La'}) then we again get a sheaf quantization similarly (Figure~\ref{fig:SLa'}):
\begin{equation}
\begin{split}
    S_{L_a'}&:=p_!\bK_{D_a'}\\
    D_a'&:=\lc (x, t)\relmid t\geq 0, x\in (a, 1-a)\rc\cup \lc (x, t)\relmid x\in (0, a], t\geq \int_x^a\lb \sqrt{a^2-(x-a)^2}-a\rb dx\rc\\
    &\cup \lc (x, t)\relmid x\in [1-a, 1), t\geq \int_{1-a}^x\lb \sqrt{a^2-(x-1+a)^2}-a\rb dx\rc
\end{split}
\end{equation}
Then we get $S_{L_a'}\rightarrow \bK_{(0,1)}\boxtimes \bK_{t\geq 0}$ as $a\rightarrow 0$. 

This suggests an interpretation of the puzzle $\SS(\bK_{(0,1)})=\SS(\bK_{[0,1]})$. A correct way to enlarge $\frakL(X, \omega)$ to include singular Lagrangians is not just by including singular Lagrangian and try to consider brane data on it. Rather, we should consider the completion of the set of Lagrangian branes with respect to the metric. Then $\SS(\bK_{(0,1)})$ and $\SS(\bK_{[0,1]})$ will correspond to different points in the completed space. In the terminology of Guillermou--Viterbo~\cite{GuillermouViterbo}, they have the same $\gamma$-support but are not the same. A relevant categorical machinery is also suggested by Fukaya~\cite{Fukayacompletion} recently.
\end{exa}

\begin{rem}
Asano--Ike~\cite{AIcomplete} and Guillermou--Viterbo~\cite{GuillermouViterbo} observed the completeness of the interleaving distance. It might be useful to think about sheaf quantization of coisotropic branes.
\end{rem}

\subsection{Brane Structure}\label{section:brane}
An object of a Fukaya category is a structure on a Lagrangian submanifold called a Lagrangian brane. We will describe it briefly, see, for examples, Seidel~\cite{Seidel} for details.

For a symplectic vector space, the set of Lagrangian subspaces form a manifold, called the Lagrangian Grassmannian manifold. It is known that the $H^1$ of the Lagrangian Grassmannian is $\bZ$, and its generator is called the Maslov class.
 
Let $X$ be a symplectic manifold. The tangent fibers of $X$ are symplectic vector spaces. Then the union of the Lagrangian Grassmannian of the tangent fibers form a fiber bundle $\mathrm{LG}_X$ over $X$. Take a class of $H^1(\mathrm{LG}_X,\bZ)$ such that restriction to each fiber induces the Maslov class. The associated covering $\widetilde{\mathrm{LG}}_X$ of $\mathrm{LG}_X$ is called the Maslov covering.

If $L\subset X$ is a Lagrangian submanifold, then the restriction of $\mathrm{LG}_X$ to $L$ has a canonical section. A lifting of this section to $\widetilde{\mathrm{LG}}_X|_L$ is called a grading of the Lagrangian submanifold. This is a part of the brane data.

Another part of the brane data is given by a relative Pin structure. For a manifold $M$, the classification map of its tangent bundle is a map
\begin{equation}
    M\rightarrow \mathrm{BO}(n).
\end{equation}
Then the second Stiefel--Whitney class is given by the composition
\begin{equation}
    M\rightarrow \mathrm{BO}(n)\rightarrow \mathrm{B}^2\bZ/2,
\end{equation}
where the second arrow is the universal second Stiefel--Whitney class.
If this is null-homotopic, since we have a homotopy exact sequence
\begin{equation}
    \mathrm{BPin^+}(n)\rightarrow \mathrm{BO}(n)\rightarrow \mathrm{B}^2\bZ/2,
\end{equation}
we have a lift of the map to $M\rightarrow \mathrm{BPin^+}(n)$. This is a Pin structure of $M$. 

For a Lagrangian submanifold, we can similarly formulate a relative Pin structure. Take an element $w\in H^2(X; \bZ/2)$. Consider the sum of the second Stiefel--Whitney class of the tangent bundle $w_2(TL)\colon L\rightarrow \mathrm{B}^2\bZ/2$ and $-w|_{L}$. A null homotopy of the sum is a relative Pin structure.
\begin{dfn}
A quadruple $(L, \cE, \alpha, s)$ is called a Lagrangian brane where $L$ is a Lagrangian submanifold, $\cE$ is a (derived) local system on $L$, $\alpha$ is a grading of $L$, $s$ is a relative Pin structure.
\end{dfn}

\begin{thm}[\cite{Guillermou, JT}]
There is a bijection between the set of sheaf quantizations of the Lagrangian submanifold $L$ and the set of Lagrangian branes.
\end{thm}

\begin{rem}
The above statement is slightly imprecise: There exists a nontrivial interaction between $\cE$ and $s$. 
\end{rem}

Another way to say the above theorem is as follows: Fix an exact Lagrangian, a grading and a relative pin structure of it. Then there exists a bijection between the set of sheaf quantizations of $L$ and the category of complexes of local systems on $L$. Sometimes one uses the term ``a sheaf quantization of $L$" only for the ones corresponding to the rank 1 local systems under the correspondence.

For a given sheaf quantization $\cE$ of a Lagrangian submanifold, we have $\SS(\cE)\bs T^*_{M\times \bR}M\times \bR=L_f$ for a primitive $f$. Reflecting $(x,\xi)\in L$, there exists a non-vanishing local cohomology of $\cE$ at $(x, f(x, \xi))$ in the codirection $(\tau\cdot \xi, \tau)$ by the definition of microsupport. The local cohomology is a priori only unique up to shift, but one can fix the shift by using the grading. Moreover, the local cohomology at different points can be glued up by using the relative pin structure. Then we get a local system over $L$. This local system is sometimes called microlocal monodromy~\cite{STZ}.

\begin{rem}
See \cite{ JT, Jin, ShendeH-principle}, etc.\  for more general treatments of brane structures and their relation to $J$-homomorphisms.
\end{rem}

\subsection{Origin of Tamarkin's Trick: Deformation Quantization}
The origin of Tamarkin's trick lies in deformation quantization, which we will attempt to explain here. The relationship between deformation quantization and sheaf quantization will be also explained in the context of $\hbar$-RH correspondence later.

The relationship between the Planck constant $\hbar$ and the extra variable was originally pointed out in the context of deformation quantization. First, it was pointed out by Kashiwara that micro-differential operators can be defined on general contact manifolds~\cite{Kashiwaracontact} by gluing. Then Polesello--Schapira~\cite{PS} used the result to construct ``micro-differential operators" on symplectic manifolds by taking the contactification of them with an extra direction $\bC_t$. However, due to this additional direction $\bC_t$, the element $\partial_t$ survives in the resulting ring, which we think as $\hbar$. D'Agnolo--Schapira~\cite{D'AgonoloSchapiraLagrangian} considered Lagrangian quantization in this context.

Based on these results, Tamarkin's lift can be interpreted as follows. For the simplicity, we consider the 1-dimensional case. Let $p,q$ be canonical coordinates of $T^*\bR$; The Poisson bracket is $\{p,q\}=1$. We denote the corresponding operators by $\widehat{q}, \widehat{p}$. Then, under the canonical quantization, $[\widehat{p}, \widehat{q}]=\hbar$. Let us represent it as $\widehat{p}=\hbar\partial=\hbar \widehat{\xi}$ where $\xi$ is the symbol of $\partial$. Furthermore, we identify $\hbar=\tau^{-1}$. Then $(q,p)\in L$ can be interpreted as $(q, \hbar \xi)=(q, \xi/\tau)\in L$ in the classical limit. This is precisely one of the conditions defining $L_f$. Another condition $\tau>0$ is making this well-defined. The third condition $t=-f(x, \xi/\tau)$ is posed to get a Lagrangian submanifold.

\begin{rem}
In K. Saito's primitive form theory~\cite{KSaito, SaitoTakahashi}, we have a counterpart of $\hbar$ as $\delta$. The definition of Gauss--Manin conneciton is similar to Tamarkin's trick.
\end{rem}

\section{Sheaf quantization version 2}\label{sectoin:v2}
\subsection{Novikov rings}
Tamarkin's trick can only handle exact Lagrangian submanifolds. A method for dealing with general non-exact Lagrangians was proposed in \cite{kuwagaki2020sheaf}~\footnote{A method for dealing with rational Lagrangians is used in \cite{AI} which replace $\bR$ by $S^1$, whose period  depends on the given Lagrangian}. The method seems to be the correct method in the sense that the Novikov ring comes out naturally.

What is the Novikov ring? It was introduced to solve the following problem: The counting of pseudo-holomorphic maps appeared in the definition of Fukaya category is not finite in general. In other words, there can be situations where the expression (\ref{eqn:2.1}) does not make sense as it is. Therefore, we make the following modifications to the formula. First, consider the nonnegative real numbers $\bR_{\geq 0}$ as a semigroup by the addition. Let $\bK$ be a field and $\bK[\bR_{\geq 0}]$ be the polynomial ring of the  semigroup $\bR_{\geq 0}$. For $a\in \bR_{\geq 0}$, we write $T^a\in \bK[\bR_{\geq 0}]$ for the corresponding indeterminate. Let $T^a\bK[\bR_{\geq 0}]$ be the ideals of $\bK[\bR_{\geq 0}]$ generated by $T^a$. The Novikov ring is defined as
\begin{equation}
    \Lambda_0:=\lim_{\longleftarrow \atop a\rightarrow \infty}\bK[\bR_{\geq 0}]/T^a\bK[\bR_{\geq 0}].
\end{equation}
An element of $\Lambda_0$ can be represented by $\sum_{a\in \bR_{\geq 0}}c_aT^a$ where $c_a\in \bK$ and the set $\lc a\in \bR_{\geq 0}\relmid c_a\neq 0\rc$ has no accumulation points.

Let us define the hom-spaces of Fukaya category as a free module over $\Lambda_0$ generated by the intersection points. We replace the equation (\ref{eqn:2.1}) with
\begin{equation}
    p_2\circ p_1=\sum_{p_3\in L_1\cap L_3}\sum_{f\in \cM(L_1, L_2, L_3, p_1, p_2,p_3)}T^{\int_{D}f^*\omega}\cdot p_3
\end{equation}
It is known that this equation is well-defined over the Novikov ring. This phenomenon is called the Gromov compactness. 

We can summarize the above as follows: Even if we start from any preferred ring $\bK$, our Fukaya category must be defined over the Novikov ring over $\bK$.

\begin{rem}
If one works with exact symplectic manifolds and exact Lagrangian submanifolds, you can work over $\bC$ due to the absence of bubbling. This is the reason why Nadler--Zaslow and Ganatra--Pardon--Shende equivalence works well over $\bC$. If one is not interested in nonexact case and quantitative phenomena, the Novikov ring/field is not very useful.
\end{rem}

\subsection{Equivariant sheaves}
The method of using equivariant sheaves for sheaf quantization was introduced in \cite{kuwagaki2020sheaf} to use Tamarkin's trick for non-exact Lagrangian while still keeping the size of the data reasonable. Let me explain this a little. Let $L$ be a general Lagrangian submanifold in $T^*M$. For a point $x$ in $L$, take a contractible neighborhood $U$ in $L$. Then the restriction of the standard Liouville form $\lambda|_U$ has a primitive function $f$. Therefore, $U_f$ can be defined in the same way as $L_f$ above.

For $c\in\bR$, the map $T_c\colon T^*M\times T^*\bR_t\rightarrow T^*M\times T^*\bR_t$ is the translation mapping
\begin{equation}
    (x, \xi, t, \tau)\mapsto (x, \xi, t+c, \tau).
\end{equation}
Then $\widetilde{U}:=\bigcup_{c\in\bR}T_cU_{f+c}$ is independent of how $f$ is taken. So, we take a contractible covering $\lc U_i\rc$ of $L$, and set
\begin{equation}
    \widetilde{L}:=\bigcup_i\widetilde{U_i},
\end{equation}
which only depends on $L$, is independent of the choice of the covering.

As with $L_f$, since $\widetilde{L}$ is a conic subset of $T^*M\times T^*\bR_t$, we can consider a sheaf with microsupport in it. Furthermore, it makes sense even if $L$ is not exact. (If $L$ is exact, it is simply $\bigcup_{c\in \bR}L_{f+c}$). However, $\widetilde{L}$ is too large and not a Lagrangian. Therefore, we reduce the information by considering $\bR$-equivariant sheaves.

For generalities of equivariant sheaves, one can refer to \cite{Tohoku, BernsteinLunts}. Consider $\bR$ as a group with respect to the usual addition. We consider it as a discrete topological group. Then, we let $\bR$ act on $M\times \bR_t$ by addition on the right factor. This is a continuous action, and we consider equivariant sheaves for this action. We denote the (derived) category of equivariant sheaves by $\Sh^\bR(M\times \bR_t)$. Furthermore, consider here also the quotient by the sheaves microsupported on $\{\tau\leq 0\}$. The resulting category is denoted by $\Sh_{\tau>0}^\bR(M\times \bR_t)$.

\begin{dfn}[Sheaf quantization version  2]
Let $\cE$ be an object of $\Sh^\bR_{\tau>0}(M\times \bR_t)$. We say it is a sheaf quantization version 2 (SQv2 for short) of $\mu\mathrm{supp}(\cE)=\rho(\SS(\cE)\cap \{\tau>0\})$. Again, we sometimes impose the microstalk conditions.
\end{dfn}
Again, the category is equipped with a Lagrangian category structure by $\mu\mathrm{supp}$.

We would like to describe the case when $M$ is a point $\{*\}$. Then the cotangent bundle is also a point, which is a 0-dimensional symplectic manifold. Inside it, we have a Lagrangian submanifold which is the same as the whole space. Then a sheaf quantization of Lagrangian submanifold will live in $\Sh_{\tau>0}^\bR(\bR_t)$. An object represented by  $\bigoplus_{c\in \bR}\bK_{t\geq c}$ is a sheaf quantization, whose equivariant structure is specified by canonical isomorphisms
\begin{equation}
    T_{c'}\bigoplus_{c\in \bR}\bK_{t\geq c}\cong \bigoplus_{c\in \bR}\bK_{t\geq c+c'}\cong \bigoplus_{c\in \bR}\bK_{t\geq c}.
\end{equation}
Then, one can calculate the endomorphism ring of the object as
\begin{equation}
    \End(\bigoplus_{c\in \bR}\bK_{t\geq c})\cong \Lambda_0.
\end{equation}
Here the isomorphism can be described as follows: As we noted earlier, we have $T^{c'}\colon \bK_{t\geq c}\rightarrow \bK_{t\geq c+c'}$ for $c'\geq 0$. Taking the direct sum, we have a morphism $T^{c'}\in \End(\bigoplus_{c\in \bR}\bK_{t\geq c})$. Under the above isomorphism, this morphism corresponds to the indeterminate $T^{c'}\in \Lambda_0$ as the notation suggests.

\begin{exa}
\begin{enumerate}
    \item Let $L$ be an exact Lagrangian in $T^*M$ and $\cE$ be an SQv1 of $L$. Then $\bigoplus_{c\in \bR}T_c\cE$ with an obvious equivariant structure gives an SQv2 of $L$.
    \item For non-exact Lagrangian submanifolds, no general existence theorem is known. Some examples of sheaf quantizations associated to deformation quantizations of Lagrangian submanifolds are constructed in \cite{kuwagaki2020sheaf, kuwagakihRH}.
\end{enumerate}
\end{exa}
Upon the above definitions, it is natural to ask a generalization of SQv2 to the case of Weinstein manifolds. It is the subject of our ongoing work with Ike~\cite{IkeKuwagaki}. Namely, we will define the category of SQv2 for Weinstein manifolds in \cite{IkeKuwagaki} and study quantitative results.

Then one can expect a version of Nadler--Zaslow--Ganatra--Pardon--Shende equivalence over the Novikov ring: There should be an appropriate definition of infinitesimally wrapped Fukaya categories of Weinstein manifolds defined over the Novikov ring, which also includes non-exact Lagrangian branes, which is equivalent to our category of sheaf quantizations. This conjectural equivalence will parametrize sheaf quantizations by the set of (unobstructed) branes, in particular, a Lagrangian category equivalence.

The equivalence should also hold in more generality: On the one hand, we can define SQv2 whose underlying Lagrangians are singular/coisotropic. We should enlarge the corresponding Fukaya category too, by, for example, Fukaya's method~\cite{Fukayacompletion}. 

Also, obstructed Lagrangian branes can define an object over $\Lambda_0/T^c\Lambda_0$ for some $c$. Such an object can be treated with Guillermou--Asano--Ike's additional variable~\cite{Guillermoucollection, AI}. We should also expect an extension of the equivalence which equates these obstructed objects.

Completing this program will give an explanation to the question of why sheaf theory can recover the qualitative/quantitative Floer-theoretic results.

\begin{rem}
There exists a variant of the construction of $\Lambda_0$ here. Again, we consider the sheaves over $\bR_t$ with positive microsupport. Instead of considering $\bR_t$-equivariance, we use $\bZ$-equivariance. We denote the resulting category by $\Sh^{\bZ}_{\tau>0}(\bR_t)$. Then the endomorphism ring $\bigoplus_{c\in \bZ}\bK_{t\geq c}$ will give a formal power series ring.

We can consider more variants by considering equivariancy with respect to subgroups of $\bR_t$.

Let us fix a prime number $p$. We set $\bK:=\bF_p$. Consider the subgroup $\bZ[1/p^\infty]$ of $\bR$ generated by $\bZ$ and $1/p^n$ for all $n$.  Let us consider equivariant sheaves withe respect to the action of $\bZ[1/p^\infty]$ on $\bR$. Then the resulting category $\Sh^{\bZ[1/p^\infty]}_{\tau>0}(\bR_t)$ contains an object represented by $\bigoplus_{c\in \bZ[1/p^\infty]}\bK_{t\geq c}$. The endomorphism ring of this object is isomorphic to $\bF_p[[t]][t^{1/p^{\infty}}]$. The quotient field is $\bF_p[[t]](t^{1/p^{\infty}})$, which often appears in number theory, in particular perfectoid theory~\cite{Scholze}.
\end{rem}

\subsection{Cone version}
Novikov rings having higher-dimensional exponents can also be useful. This is analogous to what was done in Guillermou--Schapira~\cite{GuiS} (based on Kashiwara--Schapira~\cite{KashiwaraSchapira}).

Let $\gamma$ be a strictly convex closed cone in $\bR^n$ i.e., it is closed, convex, stable under $\bR_{>0}$-scaling, and strict $\gamma\cap (-\gamma)=\{0\}$. We can consider $\gamma$ as a semigroup by the addition and denote the associated polynomial ring by $\bK[\gamma]$. For $c\in \gamma$, we denote the corresponding indeterminate by $T^c$. 
In this situation, we define the associated Novikov ring as follows. For $r>0$, we denote the ideal generated by $T^c$'s with $|c|\geq r$ by $I_r$ where $|c|$ is the usual Euclidean norm. We the set
\begin{equation}
    \Lambda^\gamma_0:=\lim_{r\rightarrow \infty}\bK[\gamma]/I_r.
\end{equation}
When $n=1$ and $\gamma=\bR_{\geq 0}$, we recover the usual Novikov ring.

Let $\gamma^\vee$ be the polar dual cone: $\gamma^\vee:=\lc x\in \bR^n\relmid \la x, y\ra\geq 0 \text{ for any $y\in \gamma$}\rc$. Consider the category of sheaves on $\bR^n$ quotiented by the sheaves microsupported on the complement of $\gamma^\vee$. We denote it by $\Sh_{\gamma>0}(\bR^n)$. Again, we can consider the $\bR^n$-equivariant version and denote it by $\Sh^{\bR^n}_{\gamma>0}(\bR^n)$. 
In this category, we have an object $\bigoplus_{c\in \bR^n}\bK_{\gamma+c}$. The endomorphism ring is isomorphic to $\Lambda_0^\gamma$.

This higher-dimensional generalization is used in the context of $\hbar$-RH correspondence~\cite{kuwagakihRH}. See the next subsection.

\begin{comment}

\subsection{Asano--Ike--Guillermou's more additional variable}
It can be observed that the sheaf quantization of the obstructed Lagrangian is not possible. However, looking at Guillermou's original construction, even if it is obstructed, it is possible to create some sheaf quantization by making it double. If we treat it as it is, it does not behave well, but if we add another dimension variable and extend it in that direction, we can obtain a good energy evaluation. In the language of Novikov rings, this is considered to be an object defined by modulo its energy.

\end{comment}

\subsection{WKB state}
Let $X$ be a symplectic manifold. Let us take a prequantum line bundle $\cL$ and fix a polarization. Then we get the state space $\cH$ of the geometric quantization (see \S \ref{section:quantization}). 

In the context of geometric quantization, a Lagrangian submanifold (+some data) should give a state (an element of the state space) via the WKB method (see, for example, ~\cite{BatesWeinstein}. Also, Tygan's \cite{Tsygan} is closer to our context).

Let us consider geometric quantization in the case when $X$ is a cotangent bundle $T^*M$. Since it is an exact symplectic manifold, one can take the trivial line bundle as a prequantum bundle. For a graph Lagrangian submanifold $L=\mathrm{Graph}(df)$,  an associated WKB state is of the form $e^{f/\hbar}\sum_{i}\psi_i\hbar^i$.
Let us concentrate on the semi-classical part $e^{f/\hbar}$. It values in $\cL\bs 0$. We consider its fiber as $\bC^\times\cong  \bC/2\pi\sqrt{-1}\bZ$. Since our geometric quantization is now trivial, we can take the lift $\widetilde{\cL}:=T^*M\times \bC$, which covers $\cL\bs 0$ by the exponential map. Then $f$ (or $f/\hbar$) takes value in $\widetilde{\cL}$. 

In this setup, we can consider our extra $\bR_t$ is a real analogue of the covering of the prequantum bundle, since our $\bR_t$ is the codomain of $f$. In other words, a sheaf quantization is a geometric/topological analogue of WKB-states of a Lagrangian.

This is more manifest in the context of $\hbar$-Riemann--Hilbert correspondence. Since the statement of the correspondence is a little complicated, we will only see some key points here. The usual Riemann--Hilbert correspondence associates a monodromy sheaf to a differential equation. In $\hbar$-Riemann--Hilbert correspondence, the assignment is enhanced to encode the asymptotics of solutions. Let us see in an example:
\begin{exa}[Schr\"odinger equation]
This example is treated in \cite{kuwagaki2020sheaf}. Let us consider a (stationary) Schr\"odingere equation $(\hbar^2\partial^2-Q)\psi=0$ on a Riemann surface $C$. The differential equation can be considered as a deformation quantization of the ``spectral curve" $\lc \xi^2-Q=0\rc\subset T^*C$. In \cite{kuwagaki2020sheaf}, we constructed a sheaf quantization of the spectral curve, using exact WKB analysis~\cite{Voros, KawaiTakei}. 

We can construct a formal WKB-solution of the form $e^{\int\frac{\sqrt{Q}}{\hbar}}\psi_i\hbar^i$ of the equation. It is a (local) WKB-state for the spectral curve. A corresponding local model for SQ is given by 
\begin{equation}
    \bigoplus_{\pm}\bigoplus_{c\in \bR}T_c\bK_{t\geq \mp\Re\int_{x_0}^x\frac{\sqrt{Q}}{\hbar}}dx,
\end{equation}
which we consider as a sheaf version of ``WKB-state". 

Exact WKB analysis tells us that the formal solutions can be upgraded to actual local solutions. One can observe that the transformation matrices of the actual solutions can be used to glue up the above local models of SQ. Then we get a global SQv2 of the spectral Lagrangian.
\end{exa}

In the above example, we construct an SQv2 from an $\hbar$-differential equation. This is quite similar to the Riemann--Hilbert correspondence, where an SQv0 (a constructible sheaf) is obtained from a differential equation. So, it is natural to expect a categorical equivalence elaborating the above example.

The assignment cannot be an equivalence of categories as it is, since the sheaf side is enriched over the Novikov ring, but the differential equation side is not. A way to resolve this discrepancy is as follows: Consider $\bC_\hbar$ be the complex plane with the standard coordinate $\hbar$. Fix a sector $S=\lc a\in \bC_\hbar\relmid |a|<r, \theta_1<\arg a<\theta_2 \rc$ for some $r, \theta_1, \theta_2$. We restrict the domain of $\hbar$ to $S$. Assume that the cone $\Cone(S)=\lc a\in \bC_\hbar\relmid \theta_1<\arg a<\theta_2 \rc$ is strict and set $\gamma:=\Cone(S)^\vee$. Then we instead use sheaf quantization with $\gamma$-microsupport as in the last subsection. On the differential equation side, we will use asymptotically well-behaved functions as coefficients. In particular, $e^{-c/\hbar}$ for $c\in \gamma$ is rapid decay over $S$. The set $e^{-c/\hbar}$ forms a ring isomorphic to $\bK[\gamma]$. Then we complete the ring to obtain the Novikov ring. This is the counterpart of the appearance of the Novikov ring. The actual choice of the coefficient ring is a little more complicated and we omit it here.

In the language of WKB analysis or resurgent analysis~\cite{Ecalle}, the exponent $c$ values in the Borel plane. Again, the extra variable is Fourier--Laplace dual of $\hbar(=\tau)$. Moreover, since $c$ is related to the area of the holomorphic disk, it seems to be compatible with the instanton interpretation of the Borel singularities.

\begin{rem}
The expected equivalence indicated in the end of \S 5.2 transforms the $\hbar$-RH correspondence to an equivalence between a category of deformation quantization modules and a Fukaya category, which should be a case of Kontsevich--Soibleman's DQ-RH correspondence in their holomorhic Floer project including \cite{KonSoibHolomorphicFloer}.
\end{rem}

\section{Appendix: Prerequisites}

\subsection{Microlocal sheaf theory}
The original source, of course, is Kashiwara--Schapira's~\cite{KashiwaraSchapira, Microlocalstudy}. If the reader understands Japanese language, we recommend an introductory survey by Ike~\cite{IkeMathlog}. Here, only a very limited overview is given.

Let us fix a base field $\bK$. Let $M$ be a manifold. We denote the set of open subsets by $\frakO_M$. Then $\frakO_M$ forms a poset by the inclusion relation. A sheaf $\cF$ in this paper is a sheaf valued in $\bK$-vector spaces (or complexes of them). That is, it is a contravariant functor $\frakO_M\rightarrow \mathrm{Vect}$ that satisfies certain local-to-global conditions. We denote the (derived) category of sheaves by $\Sh(M)$.

We would like to define several important classes of sheaves.

For a fixed $m\in \bZ_{\geq 0}$, the constant sheaf of rank $m$ is characterized by the conditions that it assigns $\bK^{\oplus m}$ for any connected open subset and the restriction map between two connected open subsets is given by the identity map $\id_{\bK^{\oplus m}}$. We denote the sheaf by $\bK_M^{\oplus m}$. If the rank is 1, we simply write it as $\bK_M$. 

A sheaf $\cF$ is a locally constant sheaf if, for any point $x\in M$, there exists a neighborhood $U$ of $x$ such that the restriction $\cF|_U$ is a constant sheaf.

A constructible sheaf is a ``piecewisely locally constant sheaf". That is, a sheaf $\cF$ is a constructible sheaf when there exists a stratification $M=\bigsqcup_{S\in \frakS}S$ of $M$ and $\cF|_{S}$ is a locally constant sheaf. Here, a stratification is a locally finite decomposition of $M$ into locally closed subsets satisfying the following condition: 
\begin{itemize}
    \item If $S, S'\in \frakS$ satisfies $\overline{S}\cap S'\neq \varnothing$, then $\overline{S}\supset S'$.
\end{itemize}
In addition, we sometimes ask that each stratum is a subanalytic set and satisfies the Whitney conditions, if necessary. 

\begin{exa}
\begin{enumerate}
    \item A locally constant sheaf is a constructible sheaf. It is constructible with respect to the trivial stratification. The only stratum is $M$ itself.
    \item For a closed submanifold $N\subset M$, we write the inclusion map by $i_N\colon N\hookrightarrow M$. Let $\cF$ be a constructible sheaf on $N$ with respect to a stratification $\cS$. 
    In this situation, we can define the push-forward $i_{N*}\cF$, which is a sheaf on $M$. For an open subset $U$ of $M$, its value is defined by $i_{N*}\cF(U)=\cF(i_N^{-1}(U))$.
    Then $i_{N*}\cF$ is a constructible sheaf with respect to a stratification $\cS\cup \{M\bs N\}$.
    \item Let us consider the case when $M=\bR$. Consider an open interval $(a,b)$ for $a< b\in \bR$. We again consider the push-forward $i_*\bK_{(a,b)}$ of the constant sheaf $\bK_{(a,b)}$ along the inclusion $i\colon (a,b)\hookrightarrow \bR$. If we take $\bR= (-\infty, a)\sqcup \{a\}\sqcup (a, b)\sqcup \{b\}\sqcup (b, \infty)$ as a stratification of $\bR$, the restriction of $i_*\bK_{(a,b)}$ to each stratum is given by
    \begin{equation}
    \begin{split}
        &(i_*\bK_{(a,b)})|_{(-\infty, a)}=\bK^{\oplus 0}_{(-\infty, a)}, (i_*\bK_{(a,b)})|_{\{a\}}=\bK_{\{a\}},(i_*\bK_{(a,b)})|_{(a,b)}=\bK_{(a,b)},\\ &(i_*\bK_{(a,b)})|_{\{b\}}=\bK_{\{b\}}, 
        (i_*\bK_{(a,b)})|_{(b, \infty)}=\bK^{\oplus 0}_{(b, \infty)}, 
    \end{split}
    \end{equation}
    Thus, it is a constructible sheaf.
    \item As a variant of the above example, consider the 0-extension (or !-extension) of $\bK_{(a,b)}$ along $i$. That is, the sheafification of the presheaf
    \begin{equation}
        U\mapsto
        \begin{cases}
        &\bK_{(a, b)}(U\cap (a,b)) \text{ if $U\subset (a,b)$}\\
        &0\text{ otherwise}
        \end{cases}
    \end{equation}
    for an open set $U$ of $\bR$. We can take the same stratification as the above example:
    \begin{equation}
    \begin{split}
        &(i_! \bK_{(a,b)})|_{(-\infty, a)}=\bK^{\oplus 0}_{(-\infty, a)}, (i_! \bK_{(a,b)})|_{a}=\bK^{\oplus 0}_{\{a\}},(i_! \bK_{(a,b)})|_{(a,b)}=\bK_{(a,b)},\\
        &(i_! \bK_{(a,b)})|_{\{b\}}=\bK^{\oplus 0}_{\{b\}}, 
        (i_! \bK_{(a,b)})|_{(b, \infty)}=\bK^{\oplus 0}_{(b, \infty)}, 
    \end{split}
    \end{equation}
    Thus, it is a constructible sheaf.

\end{enumerate}
\end{exa}

Let us see the last two examples in more detail. There is a canonical morphism from the stalk of a low-dimensional stratum to the stalk of a high-dimensional stratum. For example, if we look at the point $a\in \bR$ for $i_*\bK_{(a,b)}$, we have two morphisms:
\begin{equation}
    \bK_{\{a\}}\xrightarrow{\cong} (\bK_{(a,b)})_c, \bK_{\{a\}}\xrightarrow{0} (\bK_{(-\infty, a)}^{0})_{c'},
\end{equation}
where $c\in (a,b), c'\in (-\infty, a)$. That is, from $a$ to the left, the sheaf is ``torn off" and to the right, the sheaf is ``connected".
For $i_!\bK_{(a,b)}$, conversely, the sheaf is connected in the left direction from $a$ and torn off in the right direction. The microsupport is a way to geometrically record such a situation. That is, the microsupport $\SS(i_*\bK_{(a,b)})$ of $i_*\bK_{(a,b)}$ is the union of
\begin{enumerate}
    \item the support $[a,b]$ (which is viewed as a subset of $T^*\bR$, thinking of $[a,b]$ as being in the zero-section) of the sheaf $i_*\bK_{(a,b)}$,
    \item the positive half of the cotangent fiber over $a$, ``because" it is not connected from $a$ to the left direction,\footnote{Maybe, you think the orientation should be opposite intuitively, but this one is correct. To get the correct intuition, we should consider propagations of section, as in the formal definition of microsupport below.} and
    \item the negative half of the cotangent fiber over $b$, ``because" it is not connected from $b$ to the right direction.
\end{enumerate}
Similarly, $\SS(i_! \bK_{(a,b)})$ is computed as $-\SS(i_*\bK_{(a,b)})$. Here, the minus sign means multiplying the cotangent fiber by minus one.

In general, microsupport of a constructible sheaf is a subset of the union of conormal bundles of strata in some stratification of the sheaf. So, microsupport is a kind of information that refines the data of stratification. Historically, it was introduced as the counterpart of characteristic varieties of $\cD$-modules under the Riemann--Hilbert correspondence. But,
the notion of microsupport can be defined for any sheaves, even in the case irrelevant to differential equations.
\begin{dfn}[Microsupport]\label{microsupport}Let $\cE$ be a sheaf on $M$. Then the microsupport $\SS(\cE)$ of $\cE$ is characterized by the following: $(x, \xi)\not\in \SS(\cE)$ if there exists an open neighborhood $U$ such that, for any point $x_0\in M$, and any differentiable function $f$ with $df(x_0)\in U$, the local cohomology $\bR\Gamma_{\lc y\relmid f(y)\geq f(x_0)\rc}(\cE)_{x_0}\simeq 0$.
\end{dfn}

\begin{exa}
\begin{enumerate}
    \item For a locally constant sheaf $\cF$, the microsupport of $\cF$ is the zero-section.
    \item For a closed submanifold $i\colon N\hookrightarrow M$. The microsupport of $i_*\bK_N$ is the conormal bundle $T^*_NM$ of $N$.
\item
If $M$ is a real manifold, the microsupport of the sheaf of $C^\infty$-functions $C^\infty_M$ is the whole $T^*M$. If $M$ is a complex manifold, the microsupport of the structure sheaf $\cO_M$ is the whole $T^*M$.
\end{enumerate}
\end{exa}

\subsection{Symplectic Geometry and Quantization}\label{section:quantization}
Let $X$ be a real/complex manifold. We denote the structure sheaf by $\cO_X$, which is the sheaf of $C^\infty$-functions if $X$ is real and the sheaf of holomorphic functions if $X$ is complex.
\begin{dfn}
Let $\omega$ be a 2-form on $X$. A pair $(X,\omega)$ is a symplectic manifold when it satisfies the following two conditions.
\begin{enumerate}
\item $\omega$ is closed, and
    \item $\omega\colon TX\otimes TX\rightarrow \cO_X$ is non-degenerate
\end{enumerate}
When these conditions are satisfied, $\omega$ is called a symplectic form.
\end{dfn}
An isomorphism $X\rightarrow Y$ is a symplectomorphism if it preserves a symplectic form. A symplectic form $\omega$ induces an isomorphism $TX\cong T^*X$. Along this isomorphism, the 2-form $\omega$ induces a bivector field $P$. For $f, g\in \cO_X$, we set
\begin{equation}
    \{f, g\}=P(df\wedge dg),
\end{equation}
which is called the Poisson bracket.

\begin{exa}[Cotangent bundle]
Let $M$ be an $n$-dimensional manifold. Then the cotangent bundle $T^*M$ is equipped with the standard symplectic structure as follows: Let $\{x_i\}_{i=1}^n$ be a local coordinate of $M$. We denote the coordinate function of $T^*M$ by $\xi_i$ corresponding to $dx_i$. Then one can check that the 2-form $\sum_{i=1}^nd\xi_i\wedge dx_i$ is independent of the choice of local coordinates and gives a global symplectic 2-form. 

For local functions $f, g$, the Poisson bracket is locally written as
\begin{equation}
    \lc f,g\rc:=\sum_{i=1}^n\lb \frac{\partial f}{\partial\xi_i} \frac{\partial g}{\partial x_i}-\frac{\partial g}{\partial\xi_i} \frac{\partial f}{\partial x_i}\rb.
\end{equation}
\end{exa}

\begin{dfn}
Let $(X, \omega)$ be a symplectic manifold.
Let $H$ be a function on $X$. We define a vector field $v_H$ on $X$ by the relation
\begin{equation}
    \omega(v_H, -)=dH
\end{equation}
This is called a Hamiltonian vector field associated to the Hamiltonian function $H$. The flow generated by $v_H$ is called a Hamiltonian isotopy. For each time $t$, the resulting diffeomorphism is called a time-$t$ Hamiltonian diffeomorphism. This is a symplectomorphism. When a Hamiltonian function has compact support, its Hamiltonian isotopy/diffeomorphism is also said to be compactly-supported. 
\end{dfn}

In symplectic geometry, various important concepts are expressed through Lagrangian submanifolds.
\begin{dfn}
A submanifold $L$ of $X$ is said to be coisotropic if $TL^{\perp_{\omega}}\subset TL$, where $TL^{\perp_\omega}$ is the orthogonal complement space to $TL$ with respect to $\omega$.
A submanifold $L$ is Lagrangian if it is coisotropic and $\dim L=\frac{\dim X}{2}$.
\end{dfn}
Lagrangian and coisotropic subsets that appear in microlocal sheaf theory are often singular, in which case the condition is understood to be satisfied on a dense smooth subset~\footnote{There exists a definition of coisotropicity applicable to singular subset in \cite{KashiwaraSchapira}, which is quite recently revisited and refined in \cite{GuillermouViterbo}.}.

\begin{exa}
Given a symplectomorphism $f\colon (X, \omega_X)\rightarrow (Y, \omega_Y)$ (or more weakly, preserving symplectic forms), $\mathrm{Graph}(f)\subset X\times Y$ is a Lagrangian submanifold in $(X \times Y, \omega_X\boxplus (-\omega_Y))$.
\end{exa}

Symplectic manifolds have the aspect of mathematical abstraction of phase spaces of classical mechanics. Therefore, it is natural to consider quantization of symplectic geometry. There are various formalisms to consider quantization of symplectic manifolds. Let us first look at deformation quantization.
\begin{dfn}
A deformation quantization of $(X, \omega)$ is a product structure $\star$ on $\cO_X[[\hbar]]$ such that
\begin{enumerate}
\item we can recover $\cO_X$ modulo $\hbar$, and
\item the commutator satisfies that $f\star g-g\star f=\hbar\{f, g\}+O(\hbar^2)$.
\end{enumerate}

\end{dfn}

Lagrangian submanifolds can also be quantized in this framework.
\begin{dfn} 
Let us fix a deformation quantization ring $(\cO_X[[\hbar]], \star)$ of $X$. Let $\cM$ be a module over the deformation quantization such that $\cM/\hbar\cM$ considered as a module over $\cO_X$ has a Lagrangian support $L$. In this case, we say that $\cM$ is a quantization of $L$.
\end{dfn}

\begin{exa}[The ring of differential operators]
For the simplicity, we consider algebraically. The affine space $\bC^{2n}$ has a symplectic structure by viewing it as $T^*\bC^n$. Then the ring of functions is $\bC[x_i, p_i]$ where $i=1,...,n$. 

Consider the ring of $\hbar$-differential operators $\bC[x_i, \hbar\partial_i, \hbar]$, which satisfies $[\hbar\partial_i, x_i]=\hbar$. Then this is a deformation quantization of $\bC[x_i, p_i]$ where $p_i$ is the classical limit of $\hbar\partial_i$. Specializing at $\hbar=1$, we get the usual ring of differential operators. 

If we view the ring of differential operators as a deformation quantization, a $\cD$-module is a module over the deformation quantization. The characteristic variety of $\cD$-module is in $T^*M$, which is the correct notion of ``support" for $\cD$-modules when viewed as modules over deformation quantizations. Then quantizations of Lagrangians is precisely the class of holonomic modules.
\end{exa}

Next, let us look at geometric quantization. Let $(X, \omega)$ be a symplectic real manifold and assume that the de Rham cohomology class of $\omega$ determines an element of $H^2(X, \bZ)$. Let $\cL$ be a complex line bundle with $c_1(\cL)=\omega$. Let $\cL_{S^1}$ be the associated $U(1)$ principal bundle. Furthermore, take a $U(1)$-connection $\nabla$ whose curvature 2-form is $\omega$. 

\begin{dfn}
A polarization $\cP$ is an integrable Lagrangian distribution of $TX\otimes_\bR\bC$. We set
\begin{equation}
    \cH=\lc s\in \Gamma(\cL)\relmid \nabla_vs=0 \text{ for any $v\in \cP$}\rc.
\end{equation}
This is called the state space of geometric quantization. Under good circumstances, the state space is expected to be independent of the choice of polarization.
\end{dfn}
And then, we should observe the action of observable on the state space, etc., but we omit it here.

\begin{exa}
Consider the case when $X$ is the cotangent bundle of a manifold $M$ with the standard symplectic structure. Since the symplectic structure is exact, the Chern class is zero. We can take $\cL$ as the trivial line bundle. We then equip it with the connection $d+\lambda$ where $\lambda$ is the Liouville 1-form (see the example below). Let us take the polarization of the cotangent fiber directions. Then the state space is exactly the space of functions over $M$.
\end{exa}

\begin{rem}
As mentioned in the body of the paper, the relationship between deformation quantization and geometric quantization has recently been shed in a new light called brane quantization ~\cite{gaiotto2021probing}.
\end{rem}

\subsection{Exact symplectic geometry and wrapped Fukaya category}\label{wrapped}
Let $(X, \omega)$ be a symplectic manifold. We say $(X, \omega)$ is an exact symplectic manifold if $\omega$ is exact. If $(X, \omega)$ is exact, then the volume form $\omega^{\dim X/2}$ is also exact. Hence $X$ needs to be non-compact or compact with nonempty boundary.

Given an exact symplectic manifold, we often fix a 1-form $\lambda$ such that $d\lambda=\omega$. In this situation, we say $\lambda$ is a Liouville form. In what follows, an exact symplectic manifold refers to an exact symplectic manifold with a fixed Liouville form.

Let $(X, \lambda)$ be an exact symplectic manifold. The vector field $v_\lambda$ is defined by
\begin{equation}
    \omega(v_\lambda, -)=\lambda,
\end{equation}
which is called the Liouville vector field. The following is a very important class of exact symplectic manifolds.
\begin{dfn}
Let $D$ be a compact manifold with boundary $\partial D$, $\omega$ be a symplectic form on $D$, and $\lambda$ be a Liouville form on $\omega$. We say that $(D, \lambda)$ is a Liouville domain if $v_\lambda$ is transverse to the boundary $\partial D$ and pointing to the outside of $D$. 

Let $(X,\lambda)$ be an exact symplectic manifold without boundary. We say $(X, \lambda)$ is a Liouville manifold if there exists a subset $D$ of $X$ such that 
\begin{enumerate}
    \item $(D, \lambda|_D)$ is a Liouville domain,
    \item $X=D\cup \partial D\times [0, \infty)$ where the Liouville flow of $X$ and the positive flow along $[0,\infty)$ are identified. In other words, $(X, \lambda)$ has a cylindrical end with respect to the Liouville flow. 
\end{enumerate}
\end{dfn}

\begin{rem}
Ganatra--Pardon--Shende~\cite{Liouvillesector} introduced the notion of Liouville sector that generalizes Liouville manifolds and are more appropriate for gluing arguments.
\end{rem}

\begin{exa}
Let $M$ be a compact manifold. Then $T^*M$ is an exact symplectic manifold. There exists a standard Liouville form for the standard symplectic structure. It is locally defined by $\lambda:=\sum \xi_idx_i$.

Fix a Riemannian metric $g$ on $M$. Then we can consider the unit codisk bundle $D^*M$ and the boundary $D^*M$ is the unit cosphere bundle. Then the restriction of $\lambda$ to $D^*M$ gives a Liouville domain structure and $T^*M$ is the associated Liouville manifold.
\end{exa}

Let $(X,\lambda)$ be a Liouville manifold. By the definition, there exists a submanifold $D$ which is a Liouville domain. Then $(\partial D, \lambda|_{\partial D})$ is a contact manifold i.e., $(d\lambda|_{\partial D})^n\wedge \lambda$ is nowhere zero. 
The contactomorphism class of $(\partial D, \lambda|_{\partial D})$ is independent of the choice of $D$. We call this abstract contact manifold the contact boundary $\partial X$ of $X$.

\begin{dfn}[Reeb flow]
A vector field on a contact manifold $(Y, \lambda)$ is said to be a Reeb vector field $R_\lambda$ if it satisfies $d\lambda(R_\lambda, -)=0$ and $\lambda(R_\lambda)=1$. The associated flow is called the Reeb flow.
\end{dfn}

\begin{dfn}[Symplectic stop]
A symplectic stop is a subset of $\partial X$.
\end{dfn}
A stop plays a role to ``stop" the Reeb flow.

Now we introduce the skeleton/core of a Liouville manifold. 
\begin{dfn}Let $(X, \lambda)$ be a Liouville manifold.
\begin{enumerate}
    \item The subset of $X$ which stays in a bounded region of $X$ under the Liouville flow is called the skeleton (or core) of $X$.
    \item Fix a symplectic stop $S$ of $X$. The relative skeleton (or relative core) of $(X, S)$ is the union of the skeleton and the subset which flows to $S$ under the Liouville flow.
\end{enumerate}
\end{dfn}

\begin{exa}
For a cotangent bundle $T^*M$, its skeleton is the zero-section $T^*_MM$. For a Legendrian $S\subset \partial D^*M$, the relative skeleton is $\bR_{>0}\cdot S\cup T^*_MM$, where the $\bR_{>0}$-action is the scaling of the fibers.
\end{exa}

For a Liouville manifold $X$ with a symplectic stop $S$, one can associate the wrapped Fukaya category. It can be roughly described as follows:
\begin{enumerate}
    \item Object: Exact Lagrangian branes with conic ends; i.e., at the ends it has a form that $\partial L\times [0, \infty)\subset \partial D\times [0, \infty)$ where $\partial L=L\cap \partial D$.
    \item Morphism: For an object $L_2$, let us flow $L_2$ under the Reeb flow until $L_2$ meets the stop. We denote the resulting object $L_2'$. For objects $L_1, L_2$, the hom-space between them is generated by intersections of $L_1$ and $L_2'$. More precisely, we should formulate using limits under the Reeb flow.
\end{enumerate}
The composition and higher compositions are defined by disk counting as usual. 

\begin{exa}
In the case of cotangent bundles, the zero-section and cotangent fibers are examples of conic exact Lagrangians. When $S=\varnothing$, the associated wrapped Fukaya category is generated by a cotangent fiber and is equivalent to the category of modules over the chains of the based loop space of the base manifold~\cite{Abouzaid}.
\end{exa}

There exists a class of Liouville manifolds called Weinstein manifolds, which are tightly related to Stein manifolds. We will not recall the definition here (see \cite{Eliashberg, CieliebakEliashberg}). The important property related to us is the following:
\begin{prp}
    If a Liouville manifold $(X, \lambda)$ is Weinstein, then the skeleton is a Lagrangian subset. If a symplectic stop is Legendrian, then the relative skeleton is also a Lagrangian.
\end{prp}
Ganatra--Pardon--Shende theorem~\cite{MicrolocalMorse} (i.e., Kontsevich's conjecture) asserts that the wrapped Fukaya category of $(X, S)$ can be realized as the global section category of a sheaf defined on the relative skeleton. It is realized as a certain sheaf category by Nadler--Shende~\cite{nadlershende}.

%%%%%%%%%%%%%%%%%%%%%%%%%%%%%%%%%%%%%%%%%%%%%%%%%%

\bibliographystyle{alpha}
\bibliography{brane.bib}
\end{document}